\documentclass[11pt]{article}
\usepackage{amsfonts}
\usepackage{latexsym}
\usepackage{amsmath}
\usepackage{amssymb}
\usepackage{amsthm}

\newtheorem{theorem}{Theorem}[section]
\newtheorem{lemma}[theorem]{Lemma}

\newtheorem{proposition}[theorem]{Proposition}
\newtheorem{corollary}[theorem]{Corollary}

\theoremstyle{definition}
\newtheorem{definition}[theorem]{Definition}

\theoremstyle{remark}

\newtheorem*{note*}{Note}

\makeatletter
\newcommand{\ls}{\leqslant}

\makeatother

\begin{document}

\small

\title{\bf Relative entropy of cone measures and $L_p$ centroid bodies
\footnote{Keywords: centroid bodies, floating bodies, relative entropy, $L_p$-affine surface area, $L_p$ Brunn Minkowski
theory. 2000 Mathematics Subject Classification: 52A20, 53A15 }}

\author{Grigoris\ Paouris \thanks{partially supported by an NSF grant} and Elisabeth M. Werner 
\thanks{Partially supported by an NSF grant, a FRG-NSF grant and  a BSF grant}}

\date{}

\maketitle

\begin{abstract}
\noindent
Let  $K$ be a convex body  in $\mathbb R^n$. 
We introduce  a new affine invariant, which we call $\Omega_K$,  that can be found in three different ways:
\par
as a limit of normalized $L_p$-affine surface areas, 
\par
as the relative entropy of the cone measure of $K$ and the cone measure of $K^\circ$, 
\par
as the limit of the volume difference of $K$ and $L_p$-centroid bodies.
\par
\noindent
We investigate properties of $\Omega_K$ and of related  new invariant quantities. In particular, 
we show  new affine isoperimetric inequalities and we show a ``information  inequality"  for convex bodies.

\end{abstract}

\bigskip

\section{Introduction}

\vskip 5mm
\noindent
The starting point of our investigation was the 
study of the asymptotic behavior of the volume of $L_p$ centroid bodies as $p$ tends to infinity. This study resulted in the discovery of a new 
affine invariant,   $\Omega_K$. We then showed that the quantity $\Omega_K$
is  the relative entropy of the cone measure of $K$ and the cone measure of $K^\circ$.
Cone measures have been intensively studied in recent years
(see e.g.  Barthe/Guedon/Mendelson/Naor \cite{BGMN}, Gromov/Milman
\cite{MilGr}, Naor   \cite{Na2}  and Naor/Romik \cite{NaRo} and  Schechtmann Zinn \cite{SZ})
Finally, to our 
surprise,  $\Omega_K$  appeared again  naturally in a third way, namely as a limit of normalized $L_p$-affine surface areas. 
\par
Thus,  the invariant  $\Omega_K$ introduces a novel idea -relative entropy- into the theory of convex bodies and links concepts
from classical convex geometry like $L_p$ centroid bodies and $L_p$-affine surface area with  concepts from information theory.  Such links have already been established. Guleryuz, Lutwak, Yang and Zhang \cite{GLYZ, LYZ1, LYZ2, LYZ3, LYZ4}) use $L_p$ Brunn Minkowski theory
to develop certain entropy inequalities. Also,  classical Brunn Minkowski theory is  related to information theoretic concepts (see e.g.  \cite{ABBN1, ABBN2, Ba, BBN, CT, DCT}.

\vskip 3mm
\noindent
An important affine invariant quantity in convex geometric analysis is the affine surface area, which,  for a convex body $K \in \mathbb R^n$ is defined as 
\begin{equation}\label{p=1} 
as_{1}(K) = \int_{\partial K} \kappa^{\frac{1}{n+1}}(x) d\mu(x). 
\end{equation}
\noindent 
$\kappa (x) =\kappa_{K}(x)$ is the generalized Gaussian curvature at the boundary point $x$ of $K$ and $\mu=\mu_K$ is the surface area measure on  the boundary $\partial K$.
\noindent 
Originally a basic affine invariant
from the field of affine differential geometry, it has recently
attracted increased attention(e.g. \cite{Bar,
Lu2, MW1, SW2, WY}). It  is fundamental  in the theory of valuations
(see e.g.,   \cite{A1, A2, K, LR1}), in approximation of convex bodies by
polytopes (e.g.,  \cite{Gr2, LuSchW, SW4}) and it is the subject of the
affine Plateau problem solved in ${\mathbb R}^3$ by Trudinger and
Wang \cite{TW1, Wa}.

The definition (\ref{p=1}), at least for convex  bodies  in $\mathbb{R}^2$ and $\mathbb{R}^3$ with sufficiently smooth boundary, goes back to Blaschke \cite{Bl} and was extended to arbitrary convex bodies by e.g. \cite{Lei, Lu2,  MW1,  SW2}.
\noindent
Sch\"utt and Werner showed in \cite{SW2} that the affine surface area equals 
$$as_1(K) = \lim_{\delta \rightarrow 0} c_{n} \frac{ |K| -|K_{\delta}|  }{\delta^{\frac{2}{n+1}}} .$$
\noindent 
$c_{n}$ is a constant depending only on $n$, $|K|$ denotes the $n$-dimensional volume of $K$ and $K_{\delta}$ is the {\em convex floating body} of $K$ \cite{SW2}:
the
intersection of all halfspaces $H^+$ whose defining hyperplanes
$H$ cut off a set of volume $\delta$ from $K$.

\smallskip
\noindent 
It was shown by Milman and Pajor \cite{MP} that for ``big" $\delta$ $K_{\delta }$ is homothetic, up to a constant depending on $\delta$, to the dual of the Binet ellipsoid from classical mechanics 
and consequently  $K_{\delta }^\circ $ is homothetic to the Binet ellipsoid.

\smallskip
\noindent 
Lutwak and Zhang \cite{LZ} generalized the notion of  Binet ellipsoid and introduced the {\em $L_{p}$ centroid bodies}: For a convex body  $K$  in $\mathbb R^n$ of volume $1$ and $1\ls p \ls \infty$,  the $L_{p}$ centroid body $Z_{p}(K)$ is  this convex body that has support function 
\begin{equation} \label{def:Zp}
h_{Z_{p}(K)}(\theta) = \left(\int_{K}|\langle x, \theta\rangle |^{p} dx \right)^{1/p}.
\end{equation}

\noindent Note that in \cite{LZ} a different notation and normalization was used for the centroid body. In the present paper we will follow the notation and normalization that appeared in \cite{Pa1}.

The results of this paper deal mostly with  centrally symmetric  convex bodies $K$. Symmetry is assumed mainly because the
$L_p$ centroid bodies are symmetric by definition (\ref{def:Zp}) and used to approximate the convex bodies $K$. There exists a non-symmetric definition of $L_p$ centroid bodies in \cite{Lud} (see also \cite{HabSch}).
Using this definition,  we feel the results of the paper can be carried  over to  non-symmetric convex bodies.
\smallskip

\noindent 
In Proposition \ref{prop:Zp/Kd}  we generalize the result by Milman and Pajor mentioned above and show  that the floating body $K_{\delta} $ is  - up to a universal constant - homothetic  to the centroid body $Z_{\log_{\frac{1}{\delta}}}(K)$. 
\par
$L_p$-affine surface area, an extension of affine surface area, 
was introduced by Lutwak in the ground
breaking paper \cite{Lu1}  for $p >1$ and for general $p$ by Sch\"utt and Werner \cite{SW5}. It is now at the core of the rapidly developing  $L_p$ Brunn Minkowski theory.
Contributions here include new interpretations of $L_p$-affine
surface areas \cite{MW2, SW4, SW5, WY, WY1}, 
the study of solutions of nontrivial ordinary and, respectively, partial differential
equations (see e.g. Chen \cite{CW1}, Chou and Wang \cite{CW2},
Stancu \cite{SA1, SA2}), the study of the $L_p$ Christoffel-Minkowski problem by
Hu, Ma and Shen \cite{HMS}, characterization theorems by Ludwig
and Reitzner \cite{LR1} and the study of $L_p$ affine isoperimetric inequalities by Lutwak \cite {Lu1}  and Werner and Ye \cite{WY, WY1}.
\par
From now on we will always assume that  the centroid of a
convex body $K$ in $\mathbb R^n$ is at the origin. We  write $K\in
C^2_+$,  if $K $ has $C^2$ boundary with everywhere strictly
positive Gaussian curvature $\kappa_K$. For real  $p \neq -n$, we define  the
$L_p$-affine surface area $as_{p}(K)$ of $K$ as in \cite{Lu1} ($p
>1$) and \cite{SW5} ($p <1, p \neq -n$) by
\begin{equation} \label{def:paffine}
as_{p}(K)=\int_{\partial K}\frac{\kappa_K(x)^{\frac{p}{n+p}}}
{\langle x,N_{ K}(x)\rangle ^{\frac{n(p-1)}{n+p}}} d\mu_{ K}(x) 
\end{equation}
and
\begin{equation}\label{def:infty}
as_{\pm\infty}(K)=\int_{\partial K}\frac{\kappa _K (x)}{\langle
x,N_{K} (x)\rangle ^{n}} d\mu_{K}(x), 
\end{equation}
provided the above integrals exist. $N_K(x)$ is the outer unit
normal vector at $x$ to $\partial K$, the boundary of $K$, and $\mu_K$ is the usual surface area measure on  $\partial K$.
In particular, for $p=0$
$$
as_{0}(K)=\int_{\partial K} \langle x,N_{ K}(x)\rangle
\,d\mu_{K}(x) = n|K|.
$$
For $p=1$ we get the  classical affine surface area  (\ref{p=1}) which is independent
of the position of $K$ in space.
\par
We use the $L_p$-affine surface area to define a new affine invariant in Section 3:
\begin{equation}\label{omega}
\Omega_{K} = \lim_{p\rightarrow \infty} \left(\frac{as_{p}(K)}{n |K^\circ|}\right)^{n+p}.
\end{equation}
This is a first way how $\Omega_K$ appears.
We describe properties of this new invariant.  E.g., in  Corollary \ref{cor1} we prove the following remarkable identity (\ref{id}), which is the second way how $\Omega_K$ appears:  It   shows that
the invariant $\Omega_K$ is the exponential of the relative entropy or Kullback-Leibler 
divergence $D_{KL}$ of the cone measures $cm_K$ and $cm_{K^\circ}$ of $K$ and $K^\circ$. 
\begin{equation}\label{id}
\Omega_{K}^{1/n} = \frac{|K^{\circ}|}{|K|} \exp{\bigg( - D_{KL}(N_{K}N_{K^{\circ}}^{-1}cm_{\partial K^{\circ}} \| cm_{\partial K} ) \bigg)}. 
\end{equation}
$N_K^{-1}$ is the inverse of the Gauss map. We refer to Section 3 for its  definition and that of the  relative entropy and the cone measures. See also Gromov/Milman
\cite{MilGr} and Naor   \cite{Na2}  and Naor/Romik \cite{NaRo}   for further information on cone measures.
\par
We show that the information inequality \cite{CT} for the relative entropy of the cone measures implies an ``information inequality" for convex bodies
$$ \Omega_{K} \ls \left(\frac{|K|}{|K^{\circ}|} \right)^{n} $$ 
\noindent with equality if and only if $K$ is an ellipsoid. Independently, we can derive this inequality  from properties of the $L_p$-affine surface areas.
\par
The next proposition gives a sample of some inequalities that hold for  the affine invariant $\Omega_K$, among them an  isoperimetric inequality. More can be found in Proposition \ref{prop}. 

\vskip 2mm

\noindent
{\bf Proposition} 
\noindent 
{\em Let $K$ be a  convex body with centroid at the origin.
\par
(i) \  $\Omega_{K^{\circ}} \ls \Omega_{(\widetilde{B^n_2})^{\circ}} $
\par
(ii) \  For all $p \geq 0$, 
$
\Omega_{K} \ls \left(\frac{as_{p}(K)}{n |K^\circ|}\right)^{n+p}. 
$
\par
(iii) \   $ \Omega_{K}  \ls  \left( \frac{|K|}{|K^{\circ}|}\right)^{n} $. 
\par
If $K$ is in addition in $C^2_+$, then equality holds if and only if $K$ is an ellipsoid.  
with equality holding in (i), (ii) and (iii)  if and only if  $K$ is an ellipsoid.}

\vskip 2mm
\noindent
Proposition \ref{prop:Zp/Kd}  states that the floating body $K_{\delta} $ is  - up to a universal constant - homothetic  to the centroid body $Z_{\log_{\frac{1}{\delta}}}(K)$. This,  and the geometric interpretations of $L_p$-affine surface areas in terms of variants of the floating bodies \cite{SW5, WY, WY1}, 
led us to investigate the $L_{p}$ centroid bodies also in the context of  affine surface area.  Note the similarities in bahavior of the 
floating body and the $L_p$ centroid body. Both ``approximate" $K$ as $\delta \rightarrow 0$ respectively $p \rightarrow \infty$:
If $K$  is symmetric and of volume $1$, $Z_p(K) \rightarrow K$ as $p \rightarrow \infty$.
\par
We 
found an amazing connection between the $L_p$ centroid bodies and the new invariant $\Omega_K$ which is stated in the
following theorem for convex bodies in $C^2_+$. A forthcoming paper will address general convex bodies.
\vskip 3mm
\noindent
{\bf Theorem  \ref{theorem1}}
{\em Let $K$ be a symmetric convex body in $\mathbb{R}^n$ of  volume $1$ that is in $C^2_+$. Then
\vskip 2mm
(i)\ \  $
 \lim _{p \rightarrow \infty} \frac{p}{\log p} \left(|Z_{p}^\circ(K) |-|K^\circ |\right) = \frac{n(n+1)}{2}\  |K^\circ|.$
\vskip 2mm
 \begin{eqnarray*}
&& \hskip -13mm \mbox {(ii)}\   \lim _{p \rightarrow \infty}  p \left( |Z_p^\circ(K)| - |K^\circ|-  \frac{n(n+1)}{2p} \log p \ |Z_p^\circ(K)| \right) = \\
&& \lim _{p \rightarrow \infty}  p \left( |Z_p^\circ(K)| - |K^\circ|-  \frac{n(n+1)}{2p} \log p \ |K^\circ| \right) = \\
&& \hskip 14mm - \frac{1}{2} \ \int_{S^{n-1}}  h_K(u)^{-n}\  \log \left(2^{n+1} \pi^{n-1}
h_K(u)^{n+1 } f_K(u) \right) d \sigma(u)= \\
&&\hskip 22mm  \frac{1}{2} \ \int_{\partial K} \frac{ \kappa(x)}{\langle x, N(x) \rangle^n} \   \log \left(\frac{\kappa(x)}{2^{n+1} \pi^{n-1}
\langle x, N(x) \rangle^{n+1}} \right) d\mu(x)
\end{eqnarray*}
}

\vskip 2mm
In view of Proposition \ref{prop:Zp/Kd}, the first part of the Theorem  \ref{theorem1}  came as a surprise to
us because it reveals a different behaviour of the bodies $K_{\delta}$ and $Z_{\log{\frac{1}{\delta}}}(K)$  when $\delta \rightarrow 0$. Indeed, it was shown in \cite{MW2} that 
\begin{eqnarray*}\label{as-1}
\lim _{\delta \rightarrow 0} \ c_n \frac{|(K_{\delta }
)^\circ|-|K^\circ|}{\delta ^{\frac{2}{n+1}}} =  as_{-n(n+2)}(K) =  as_{-\frac{n}{n+2}}(K^\circ) 
\end{eqnarray*}
where  $c_n$ is a constant that depends on $n$ only. 
\vskip 2mm
Even more surprising  is the second part of Theorem \ref{theorem1}.   Indeed, Proposition \ref{prop} states  that
$$
\Omega_{K} =\exp{ \left( -\frac{1}{|K^{\circ}|} \int_{\partial K} \frac{\kappa_{K}(x)}{\langle x, N_{K}(x)\rangle^{n}} \log{ \frac{ \kappa_{K}(x)}{\langle x, N_{K}(x)\rangle^{n+1}} } d\mu_{K}(x) \right) }.
$$
This,  together with Theorem \ref{theorem1} shows how the new invariant and the $L_p$ centroid bodies are related, namely for a  symmetric convex body $K$ of volume $1$ in $C^{2}_{+}$
\begin{eqnarray*}
&&\hskip -20mm \lim_{p \rightarrow \infty} \frac{2p}{n} \left( \frac {|Z_{p}^{\circ}(K)|}{|K^{\circ}|} - \left(1-\frac{n(n+1)\log{p}}{2p}\right) \right)  = \\
&&\hskip -10mm \lim_{p \rightarrow \infty} \frac{2p}{n} \left( \frac{(1-\frac{n(n+1)\log{p}}{2p})|Z_{p}^{\circ}(K)|}{|K^{\circ}|} -1 \right) = \\
&&\hskip 30mm   -\frac{1}{2} \log{\frac{\Omega^\frac{1}{n}_{K}}{2^{n+1}\pi^{n-1}} }.
\end{eqnarray*}
This is the third way how $\Omega_K$ appears.

\vskip 4mm

\noindent 
{\bf Further notation}. We work in ${\mathbb R}^n$, which is
equipped with a Euclidean structure $\langle\cdot ,\cdot\rangle $.
We denote by $\|\cdot \|_2$ the corresponding Euclidean norm, and
write $B_2^n$ for the Euclidean unit ball, and $S^{n-1}$ for the
unit sphere. Volume is denoted by $|\cdot |$. 
We write $\sigma $ for the rotationally
invariant surface measure on $S^{n-1}$. 
\par
A convex body is a compact convex subset $C$ of ${\mathbb R}^n$ with
non-empty interior. We say that $C$ is $0$ - symmetric, if $x\in
C$ implies that $ -x\in C$. We say that $C$ has centre of mass at the
origin if $\int_C\langle x,\theta\rangle dx=0$ for every $\theta\in
S^{n-1}$. The support function $h_C:{\mathbb R}^n\rightarrow
{\mathbb R}$ of $C$ is defined by $h_C(x )=\max\{\langle x,y\rangle
:y\in C\}$. 
The polar body $C^{\circ }$ of $C$ is
$C^{\circ}  = \{ y\in {\mathbb R}^n: \langle x,y\rangle \ls 1\;\hbox{for all}\; x\in C\}.
$ 
\par
Whenever we write $a\simeq b$, we mean that there
exist absolute constants $c_1,c_2>0$ such that $c_1a\ls b\ls c_2a$.
The letters $c,c^{\prime }, c_1, c_2$ etc. denote absolute positive
constants which may change from line to line. We refer to the books
\cite{Pi} and \cite{Sch} for basic facts from the
Brunn-Minkowski theory and the asymptotic theory of finite
dimensional normed spaces. 
\vskip 2mm
The authors would like to thank the American Institute of Mathematics. The idea for the paper originated 
during a stay at AIM.

\vskip 4mm

\section{Comparison of Floating bodies and $L_{p}$ centroid bodies}

\smallskip
It is well known from mechanics that the body $Z_{2}(K)$  is an ellipsoid. Its polar body $Z_{2}^{\circ}(K)$ is called the Binet ellipsoid of inertia. 
$Z_{1}(K)=Z(K)$ is the classical centroid body and it is  a zonoid by definition (see \cite{Ga, Sch}).

\noindent 
The isotropic contant $L_K$ of a convex body $K \in \mathbb R^n$ is defined as 
$$ 
L_{K} = \left(\frac{|Z_{2}(K)|}{|B^n_2|} \right)^{1/n} 
$$
\noindent 
$L_K$ is an affine invariant  and  $L_{K}\geq L_{B^n_2}$. 
\par
A major open problem in convex geometry asks if there exists a universal constant $C>0$ such that $L_{K} \leq C$.
The best -up to date- known result is due to  Klartag  \cite{Kl} and states that $L_{K}\leq C n^{\frac{1}{4}}$,  improving by a factor of logarithm an earlier result by Bourgain \cite{Bou1}.

\noindent 
Let us briefly state some of the known properties of the $L_p$ centroid bodies. For the proofs and further references see \cite{Pa1}. 
\vskip 3mm
Let $T \in SL(n)$, i.e. $T: \mathbb{R}^n \rightarrow \mathbb{R}^n$  is a linear operator with determinant $1$. Let $T^*$ denote its adjoint. Then 
$$ h_{Z_p(TK)} (\theta) =\left(\int_{TK} |\langle x, \theta \rangle|^p dx \right)^{1/p} = \left(\int_{K} |\langle x, T^{\star}(\theta) \rangle|^p dx  \right)^{1/p} =h_{Z_p(K)} (T^{\star}(\theta))$$ 
or
$$ h_{Z_p(TK)} (\theta) = h_{T(Z_p(K))} (\theta)  $$
\noindent  
By H\"older's inequality,  we have for $1\ls p\ls q\ls \infty$ that
\begin{equation}\label{inclusion1} 
Z_{1}(K) \subseteq Z_{p} (K)  \subseteq Z_{q} (K)  \subseteq Z_{\infty} (K) = K. 
\end{equation}
\noindent 
As an  application of the Brunn-Minkowski inequality, one has for $1 \ls p \ls q < \infty$ that
\begin{equation}\label{subset1}
Z_{q} (K)  \subseteq c\ \frac{q}{p}\  Z_{p} (K).
\end{equation}
 $c>0$ is a universal constant.
\par
 \noindent
Inequality (\ref{subset1}) is sharp with the right constant for the $l^1_n$-ball  \cite{FPS}.
\vskip 3mm
\noindent 
By Brunn's principle we get for $p\geq n$ and a (new)  absolute constant $c>0$ (e.g., \cite{Pa2})
\begin{equation} \label{Brunn}
Z_{p}(K) \supseteq c \  K.
\end{equation}

\vskip 4mm

\noindent 
Lutwak, Yang and Zhang \cite{LYZ} and Lutwak and Zhang  \cite{LZ} proved the following $L_p$ versions of
Blaschke Santal\'{o}   inequality and Busemann-Petty inequality.  See also Campi and Gronchi \cite{CG} for an alternative proof.
\vskip 3mm
\noindent
\begin{theorem}\cite{LYZ, LZ} \label{LYZ}
Let $K$ be a convex body in $\mathbb R^n$ of volume $1$. Then for every $1\ls p \ls \infty$ 
\begin{equation} \nonumber
 |Z_p^\circ(K)| \le |Z_{p}^\circ(\widetilde{B^n_2})|
\end{equation}
\begin{equation} \nonumber
 |Z_p(K)| \ge |Z_{p}(\widetilde{B^n_2})|
\end{equation}

with equality if and only if $K$ is an ellipsoid. 
\end{theorem}
\par
\noindent 
A computation shows that $|Z_p(\widetilde{B^n_2}|^{1/n} \simeq \sqrt{\frac{p}{n+p}} $. 
Hence,  the following inequality, proved in \cite{Pa1}  for all $p\geq 1$ and  a universal constant $c>0$,  can be viewed as an  ``Inverse Lutwak-Yang-Zhang  inequality"
\begin{equation}\label{inverseLYZ}
|Z_p(K)|^{1/n} \leq c\sqrt{\frac{p}{n+p}} L_K. 
\end{equation}

\vskip 4mm
\noindent
We now want to compare  $L_p$ centroid bodies and floating bodies. As $K$ is symmetric and has volume $1$, the floating body $K_\delta$,  for $\delta \in [0,1]$,  may be defined in the following way \cite{SW2}
 \begin{equation}\label{floating}
K_{\delta}  = \bigcap_{\theta \in S^{n-1}} \{ x\in K: |\langle x, \theta \rangle| \leq t_{\theta} \} 
\end{equation}
where $t_{\theta} = \sup \{t>0 : | \{x \in K: |\langle x, \theta \rangle| \leq t\}|=1- \delta \} $. Hence for every $\theta \in S^{n-1}$ one has that 
\begin{equation}\label{def:Kd}
h_{K_{\delta}}(\theta)= t_{\theta}. 
\end{equation}

\vskip 3mm

\begin{theorem}\label{prop:Zp/Kd}
\noindent 
Let $K$ a symmetric convex body in $\mathbb R^n$ of volume $1$. Let  $\delta \in (0,1)$. Then we have for every $\theta \in S^{n-1}$ 
$$ 
c_1 h_{Z_{\log{\frac{1}{\delta}}}(K)} (\theta) \ls h_{K_{\delta}} (\theta) \ls c_2 h_{Z_{\log{\frac{1}{\delta}}}(K)} (\theta) 
$$
or, equivalently
$$ 
c_1 Z_{\log{\frac{1}{\delta}}}(K) \subseteq K_{\delta} \subseteq c_2 Z_{\log{\frac{1}{\delta}}}(K), 
$$
\noindent 
where $c_1, c_2 >0$ are universal constants.
Consequently
\begin{equation} \nonumber
\label{eq2:prop1}
\frac{1}{c_1} Z^\circ_{\log{\frac{1}{\delta}}}(K) \supseteq K^\circ_{\delta} \supseteq \frac{1}{c_2} Z^\circ_{\log{\frac{1}{\delta}}}(K)
\end{equation}

\end{theorem}

\vskip 3mm

\noindent 
{\bf Proof.}
\par
\noindent 
For $\delta \in (c_0,1]$, $c_0$ appropriately chosen, the theorem was already shown in \cite{MP}. We assume
that $\delta \leq c_0 <1$. 
We apply Markov's inequality in (\ref{def:Zp}) and get
$$ 
| \{ x \in K : |\langle x, \theta \}| \ge e h_{Z_p(K)}(\theta) \}| \ls e^{-p}.
 $$
\noindent 
Then (\ref{def:Kd})  gives  for all $p\geq 1$,
\begin{equation}\label{est2}
e h_{Z_p(K)}(\theta) \ge h_{K_{e^{-p}}} (\theta). 
\end{equation}
\noindent 
For the other side we will use the Paley-Zygmund inequality: If $Z\ge0$ is a random variable with finite variance and $\lambda \in (0,1)$ then
$$ 
Pr\{ Z\ge \lambda E(Z)\} \ge (1-\lambda)^2 \frac{E(Z)^2}{E(Z^2)}. 
$$
\noindent 
Hence for  $Z=|\langle x, \theta\rangle|^p $ we get
\begin{equation}\label{est1} 
|\{x\in K: |\langle x, \theta \rangle|^p \ge \lambda \int_{K} |\langle x, \theta \rangle |^p dx \}| \geq (1-\lambda)^2 \frac{ \left(\int_{K} |\langle x, \theta \rangle |^p dx \right)^2 }{ \int_{K} |\langle x, \theta \rangle |^{2p} dx }.
\end{equation}
\par
\noindent
(\ref{subset1}) implies that $h_{Z_{2p}(K)}(\theta) \ls 2c h_{Z_{p}(K)}(\theta)$ , for all $\theta \in S^{n-1}$. So
$$ 
\frac{ \left(\int_{K} |\langle x, \theta \rangle |^p dx \right)^2 }{ \int_{K} |\langle x, \theta \rangle |^{2p} dx } \geq \left(\frac{1}{2c}\right)^{2p}.
$$
\par
\noindent 
Choose $\lambda =\frac{1}{2}$.  Then (\ref{est1}) becomes
$$ |\{x\in K: |\langle x, \theta \rangle| \ge \frac{1}{2} h_{Z_p(K)}(\theta) \}| \geq e^{-c_1 p}.$$
Now we use again  (\ref{def:Kd}) to get
$$ \frac{1}{2} h_{Z_p(K)}(\theta) \ls h_{K_{e^{-c_1 p}}} (\theta)  $$
or
\begin{equation}\label{est3}
h_{K_{e^{- p}}} (\theta) \ge \frac{1}{2} h_{Z_{\frac{p}{c_1}}(K)}(\theta) \ge c_2 h_{Z_{p}(K)}(\theta),
\end{equation}
where we have used  (\ref{subset1})  again.
\noindent 
(\ref{est2}) and (\ref{est3}) then imply that
$$
c_2\  h_{Z_{p}(K)}(\theta) \leq h_{K_{e^{- p}}} (\theta) \leq e\  h_{Z_{p}(K)}(\theta).
$$
Now choose $p=\log{\frac{1}{\delta}}$. This gives the theorem

\vskip 3mm
One does not expect that floating bodies and $L_{q}$ centroid bodies are identical in general. 
Indeed, observe that for $p<\infty$ the bodies $Z_{p}(K)$ are $C^{\infty}$. However one can easily check that the floating body of the cube has points of non-differentiability on the boundary. 
\smallskip

\smallskip

Theorem \ref{prop:Zp/Kd} allows us to ``pass'' results about $L_{p}$ centroid bodies to  floating bodies. In particular, (\ref{inclusion1}) and (\ref{Brunn}) imply that for $\delta< e^{-n}$, $K_{\delta}$ is isomorphic to $K$:

$$ K_{\delta}  \subseteq K  \subseteq c_{1}  K_{\delta}. $$

\noindent Moreover, (\ref{inclusion1}) and (\ref{subset1}) imply that
$$ K_{\delta_{2}} \subseteq K_{\delta_{1}}  \subseteq c_{2} \frac{\log{\frac{1}{\delta_{1}} }} { \log{\frac{1}{\delta_{2}} }} K_{\delta_{2}} \ , \  \  \mbox{for} \  \   \delta_{1} \ls \delta_{2} \ ,
$$

\noindent where $c_{1}, c_{2} >0 $ are universal constants. 

\smallskip

\noindent 
As a consequence  we get the following corollary.  There,  $d(K,L)$, resp. $ d_{BM}(K,L)$,  mean the geometric, 
resp. Banach-Mazur  distance of two convex bodies $K$ and $L$
$$
d(K,L) = \mbox{inf}\{ a \cdot b: \frac{1}{a} K \subset L \subset b K\}
$$
$$
d_{BM}(K,L) =\mbox{inf}\{ d\big(K,T(L)\big): \mbox{T is a linear operator}\}
$$
It is known that one may choose a $T \in SL(n)$ such that $T(K_{1/2})$ is isomorphic to $B^n_2$ (see \cite{MP} for details).
\begin{corollary}
\noindent Let $K$ be a symmetric convex body of volume $1$. Then for every $\delta \in (0,1) $ one has
$$  d_{BM} \left( K_{\delta}, B_{2}^{n} \right) \ls  c_{1} \log{\frac{1}{\delta}} \  ,$$
and
$$ d \left( K_{\delta},  K\right) \simeq d \left( K_{\delta},  K_{e^{-n}} \right) \ls c_{2}\frac{n}{\log{\frac{1}{\delta}}} \ , $$

\noindent where $c_{1}, c_{2} >0 $ are universal constants. 
\end{corollary}
\noindent 
Let us note that Theorem \ref{LYZ} and (\ref{inverseLYZ}) imply sharp (up to $L_K$) bounds for the volume of $K_{\delta}$. Namely, letting $c_\delta= \mbox{max}\{ \mbox{log}\frac{1}{\delta},1 \}$,

$$   c_1\  \sqrt{ \frac{c_{\delta} } {n+c_{\delta}} } \ls \left| K_{\delta}\right |^{1/n} \ls
 c_2 \  \sqrt{ \frac{c_{\delta} } {n+c_{\delta}} }\   L_K  \ ,$$

\noindent where $c_{1}, c_{2} >0 $ are universal constants. 
\medskip
\vskip 3mm
\noindent
{\bf Remark.} The corollary is also true for non symmetric $K$.
\medskip

\noindent 
In view of  a result of R. Latala and J. Wojtaszczyk \cite{LW}, Theorem \ref{prop:Zp/Kd} has another consequence: The floating body of a symmetric convex body $K$ corresponds to a level set of the Legendre transform of the logarithmic Laplace transform on $K$.
\par
\noindent
Let $x \in \mathbb R^n$ and $K$ a symmetric convex body of volume $1$. 
Let 
$$ \Lambda_{K}^{\ast}(x) := \sup_{ u\in \mathbb R^n} \bigg\{ \langle x, u \rangle - \log{\int_{K} e^{\langle x,u \rangle } d x} \bigg\} \ . $$
be the Legendre transform of the logarithmic Laplace transform on $K$.
\noindent 
For any $r>0$, let $B_{r}(K)$ be the convex body defined as
$$ B_{r}(K) := \{ x\in \mathbb R^n:  \Lambda_{K}^{\ast}(x) \ls r \}. $$

\noindent
 It was proved in \cite{LW} that $B_{p}(K)$ is isomorphic to $Z_{p}(K)$, 
$$ c_{1} Z_{p}(K) \subseteq B_{p}(K) \subseteq c_{2}Z_{p}(K), $$
\noindent where $c_{1}, c_{2} >0$ are universal constants.
\vskip 3mm
We combine  this  with Theorem \ref{prop:Zp/Kd}  and get the following 
\par
\begin{proposition}
\noindent Let $K$ a symmetric convex body of volume $1$ in $\mathbb R^n$. Then for every $\delta \in (0,\frac{1}{2} ) $ one has that 
$$ c_{1} \bigg\{ x\in \mathbb R^n:  \Lambda_{K}^{\ast}(x) \ls \log{\frac{1}{\delta}}  \bigg\}  \subseteq K_{\delta} \subseteq c_{2} \bigg\{ x\in \mathbb R^n:  \Lambda_{K}^{\ast}(x) \ls \log{\frac{1}{\delta}}  \bigg\}. $$
$c_{1}, c_{2} >0$ are universal constants.
\end{proposition}

\vskip 5mm

\section{Relative entropy of cone measures and related inequalities}

\smallskip

Let $K$ be a convex body in $\mathbb R^n$ with its centroid at  the origin. For real  $p \neq -n$, 
$L_p$-affine surface area $as_{p}(K)$ of $K$ was defined in 
(\ref{def:paffine})  and (\ref{def:infty}) in the introduction.
\par
If  $K$ is in $C^2_+$,  then (\ref{def:paffine})
and (\ref{def:infty}) can be written as integrals over the boundary
$\partial B^n_2=S^{n-1}$ of the Euclidean unit ball $B^n_2$  in $\mathbb R^n$
$$
as_{p}(K)=\int_{S^{n-1}}\frac{f_{K}(u)^{\frac{n}{n+p}}}
{h_K(u)^{\frac{n(p-1)}{n+p}}}
d\sigma(u)
$$
and
\begin{equation}\label{inf-aff}
as_{\pm\infty}(K)
=\int_{S^{n-1}}\frac{1}{h_K(u)^{n}}
d\sigma(u)
=n|K^{\circ}|.
\end{equation}
$f_K(u)$ is the curvature function, i.e. the reciprocal of the Gauss curvature $\kappa(x)$ at that point $x$ in $\partial K$ that
has $u$ as outer normal.
\vskip 3mm
We recall first  results   proved in \cite{WY}.
\begin{proposition}\cite{WY} \label{propWY}
\noindent Let $K$ be a convex body in $\mathbb {R}^n$  such that $\mu\{x \in \partial K: \kappa(x)=0 \}=0$. Let $p\neq -n $ be a real number. Then
\vskip 2mm
(i)  The function $p \rightarrow \left(\frac{as_{p}(K)}{as_{\infty}(K)}\right)^{n+p} $ is  decreasing in $p\in (-n, \infty)$.
\vskip 2mm
(ii)  The function $p \rightarrow \left(\frac{as_{p}(K)}{n |K^\circ |}\right)^{n+p} $ is   decreasing in $p\in (-n, \infty)$. 
\vskip 2mm
(iii) The function $p \rightarrow \left(\frac{as_{p}(K)}{n|K|}\right)^{\frac{n+p}{p}}$ is  increasing in $p\in (-n, \infty)$. 
\vskip 2mm
(iv) $as_{p}(K)= as_{\frac{n^{2}}{p}}(K^{\circ}) $.
\end{proposition}
\vskip 4mm
\noindent
{\bf Remarks.}
\par
(i) It was shown in \cite{Hug} that for $p>0$ (iv) holds without any assumptions on the boundary of $K$.  
\vskip 2mm
(ii) 
Also, it follows from the proof in \cite{WY} that (i), (ii) and (iii) hold without assumptions on the boundary of $K$ if $p \geq 0$.
\vskip 2mm
(iii) Proposition \ref{propWY} (ii) is not explicitly stated in \cite{WY},  but follows (without any assumptions on the boundary of $K$ if $p \geq 0$) from e.g., inequality (4.20) of  \cite{WY} and the following fact  (see \cite{SW5}): Let $K$ be a convex body in $ \mathbb{R}^n$. Then
\begin{equation}\label{equal}
as_{\infty}(K) \leq n|K^ \circ|
\end{equation} with equality if   $K$ is in $C^2_+$.
\vskip 2mm
(iv) Strict monotonicity and characterization of equality in Proposition \ref{propWY} (i), (ii) and (iiii):
\newline
Proposition \ref{propWY} (i), (ii) and (iii)  -without equality characterization- was proved in  \cite{WY} using 
H\"older's inequality. 
It follows immediately from the characterization of equality in 
H\"older's inequality, that strict monotonicity holds in \ref{propWY} (i), (ii) and (iii) if and only if $\mu$ -a.e.  on $\partial K$
\begin{equation*}
\frac{\kappa(x)} {\langle x, N(x) \rangle ^{n+1} } = c,
\end{equation*}
where $c>0$ is a constant - 
unless $\kappa(x)=0$ $\mu$ -a.e. on $\partial K$. If $\kappa(x)=0$ $\mu$ -a.e. on $\partial K$, then 
for all $p > 0$, $\left(\frac{as_{p}(K)}{as_{\infty}(K)}\right)^{n+p}=\mbox{constant}=0$, $\left(\frac{as_{p}(K)}{n |K^\circ |}\right)^{n+p}=\mbox{constant}=0$ and $\left(\frac{as_{p}(K)}{n|K|}\right)^{\frac{n+p}{p}}=\mbox{constant}=0$.
\par
If $K$ is in $C^2_+$, then the following theorem due to Petty \cite{Pe} implies that we have strict monotonicity in \ref{propWY} (i), (ii) and (iii) unless $K$ is an ellipsoid, in which case the quantities in \ref{propWY} (i), (ii) and (iii) are all constant equal to $1$.

\vskip 4mm
\begin{theorem} \cite{Pe} \label{Petty} 
Let $K$ be a convex body in $C^2_+$. $K$ is an ellipsoid if and only if  for all $x$  in $\partial K$
$$
\frac{\kappa(x)} {\langle x, N(x) \rangle ^{n+1} } = c,$$
where $c>0$ is a constant.
\end{theorem}

\vskip 4mm

We now introduce  new affine invariants.
\vskip 2mm
\noindent
\begin{definition} \label{newinv}
\par
\noindent
(i) Let $K$ a convex body in $\mathbb R^n$ with centroid at the origin. We define
$$ \Omega_{K} = \lim_{p\rightarrow \infty} \left(\frac{as_{p}(K)}{n |K^{\circ}|)}\right)^{n+p},$$
\par
\noindent
(ii) Let  $K_1, \dots, K_n$ be convex bodies in $\mathbb R^n$, all  with centroid at the origin. We define
$$ \Omega_{K_1, \dots K_n} = \lim_{p\rightarrow \infty} \left(\frac{as_{p}(K_1, \dots, K_n)}{as_{\infty}(K_1, \dots, K_n)}\right)^{n+p}.$$
Here
$$
as_{p}(K_1, \dots, K_n) = \int
_{S^{n-1}}\bigg[h_{K_1}(u)^{1-p}f_{K_1}(u)\cdots
h_{K_n}^{1-p}f_{K_n}(u)\bigg]^{\frac{1}{n+p}}\,d\sigma (u)
$$
is the mixed $p$-affine surface area introduced for $1 \leq p < \infty$ in \cite{Lu1}  and for general $p$ in  \cite{WY1}.

\begin{eqnarray*}
as_{\infty}(K_1, \dots, K_n) &=& \int _{S^{n-1}}
\frac{1}{h_{K_1}(u)} \cdots \frac{1}{h_{K_n}(u)} \,d \sigma (u) \\
&=& n \tilde{V}(K_1^\circ, \cdots, K_n^\circ)
\end{eqnarray*}  
is the dual mixed volume of $K_1^\circ, \cdots, K_n^\circ$,  introduced by Lutwak in \cite{Lut1975}.
\end{definition}

\vskip 3mm
\noindent
We will concentrate on describing the properties of $ \Omega_{K}$. The analogue properties for the invariant $ \Omega_{K_1, \dots K_n}$ also hold and are proved  similarly using results about the mixed $p$-affine surface areas proved in \cite{WY1}. For instance, the analogue to Proposition \ref{prop} (ii)
holds: For all $ p \geq 0$
$$
\Omega_{K_1, \dots K_n} \leq \left(\frac{as_{p}(K_1, \dots, K_n)}{as_{\infty}(K_1, \dots, K_n)}\right)^{n+p}.
$$
This follows from a monotonicity behavior of $\left(\frac{as_{p}(K_1, \dots, K_n)}{as_{\infty}(K_1, \dots, K_n)}\right)^{n+p}$ which was shown in \cite{WY1}.
And the  analogue to Proposition \ref{prop2} (ii)
holds
$$
\Omega_{K_1, \dots K_n} = \exp{ \left( \frac{1}{as_{\infty}(K_1, \dots, K_n)} \int_{S^{n-1}} \frac{\sum_{i=1}^n \log\left[f_{K_i} h_{K_i}^{n+1}\right]}{\prod_{i=1}^{n} h_{K_i}}d\sigma \right) }
$$
\vskip 3mm
\noindent
{\bf Remarks.}
\vskip 2mm
(i)
If $\mu\{x \in \partial K: \kappa(x)=0 \}=0$, then $\Omega_{K} > 0$. If $\kappa(x)=0$ $\mu$ -a.e.  on $\partial K$, then  $\Omega_{K} = 0$. In particular,
$\Omega_{P} = 0$ for all polytopes $P$.
\vskip 2mm
(ii) If $K$ is in $C^2_+$, then, by (\ref{equal}),  $as_{\infty}(K) = n|K^ \circ|$ and thus we then also have 
\begin{equation}\label{def2}
\Omega_{K} = \lim_{p\rightarrow \infty} \left(\frac{as_{p}(K)}{as_{\infty}(K)}\right)^{n+p}.
\end{equation}
\vskip 2mm
(ii)  As for all $p \neq -n$ and for all linear, invertible transformations $T$,  
$as_p(T(K))= |\mbox{det} (T)|^\frac{n-p}{n+p} as_p(K)$  (see \cite{SW5}) and $ as_{p}(T(K_1), \dots, T(K_n)) = |\mbox{det} (T)|^\frac{n-p}{n+p} as_{p}(K_1, \dots, K_n)$ \cite{WY1},  we get that  
\begin{equation}\label{invariance}
\Omega_{T(K)} =  |\mbox{det} (T)|^{2n} \  \Omega_K.
\end{equation}
and 
\begin{equation*}
\Omega_{(T(K_1), \dots, T(K_n))} =  |\mbox{det} (T)|^{2n} \  \Omega_{K_1, \dots K_n}.
\end{equation*}

In particular, $\Omega_K$  and $ \Omega_{K_1, \dots K_n} $ are  invariant under linear  transformations $T$ with $|\mbox{det} (T)|=1$.

\vskip 4mm
\noindent
\begin{corollary} \label{cor}
\noindent Let $K$ be a  convex body  $\mathbb R^n$ with centroid at the origin.
Then 
\vskip 3mm
\noindent
 $$ 
 \Omega_{K} =  \lim_{p\rightarrow 0}  \left(\frac{as_{p}(K^{\circ})}{n|K^{\circ}|}\right)^{\frac{n(n+p)}{p}}. 
 $$
\end{corollary}
\vskip 3mm
\noindent{\bf Proof.}
By  Proposition \ref{propWY} (iv) and Remark (i) after it
\begin{eqnarray*}
\Omega_{K} &= & \lim_{p\rightarrow \infty} \left(\frac{as_{p}(K)}{n |K^{\circ}|}\right)^{n+p} =  \lim_{p\rightarrow \infty} \left(\frac{as_{\frac{n^2}{p}}(K^{\circ})}{n|K^{\circ}|}\right)^{n+p} \\
&=&  \lim_{q\rightarrow 0} \left(\frac{as_{q}(K^{\circ})}{n|K^{\circ}|}\right)^{n+\frac{n^2}{q}} = \lim_{q\rightarrow 0} \left(\frac{as_{q}(K^{\circ})}{n|K^{\circ}|}\right)^{\frac{n(n+q)}{q}} 
\end{eqnarray*}
\vskip 3mm
\noindent
{\bf Example.}
\par
\noindent
For $1 \leq  r <\infty$, let $B^n_r=\{x \in \mathbb{R}^n: \left(\sum_{i=1}^{n} |x_i|^r \right)^\frac{1}{r} \leq 1\}$
and let $B^n_\infty=\{x \in \mathbb{R}^n: \mbox{max}_{1 \leq i \leq n} |x_i| \leq 1\}$.
Then a straightforward, but tedious calculation gives
\begin{equation}\label{lpball}
\Omega_{B^n_r} = \frac{\mbox{exp}\bigg(- \frac{n^2(r-2)}{r} \left( \frac{\Gamma^{\prime}(\frac{r-1}{r})}{\Gamma(\frac{r-1}{r})} - \frac{\Gamma^{\prime}(n\frac{r-1}{r})}{\Gamma(n\frac{r-1}{r})}\right)\bigg)}{(r-1)^{n(n-1)}}.
\end{equation}
Indeed, it was shown in  \cite{SW5} that  $as_p(B^n_r)=\frac{2^n(r-1)^\frac{p(n-1)}{n+p}}{r^{n-1}} \frac{\left(\Gamma(\frac{n+rp-p}{r(n+p)}\right)^n}{\Gamma(\frac{n(n+rp-p)}{r(n+p)})}$. 
Therefore
$$
\frac{as_{p}(B^n_r)}{n |(B^n_r)^\circ|} = \frac{1}{(r-1)^\frac{n(n-1)}{n+p}}\ 
\frac{\left(\Gamma(\frac{n+rp-p}{r(n+p)}\right)^n}{\Gamma(\frac{n(n+rp-p)}{r(n+p)})} 
\  \frac{\Gamma(\frac{n(r-1)}{r})}{\left(\Gamma(\frac{r-1}{r})\right)^n}
$$
and
$$
\Omega_{B^n_r} = \left(\frac{as_{p}(B^n_r)}{n |(B^n_r)^\circ|} \right)^{n+p}=
\frac{\mbox{exp}\bigg(- \frac{n^2(r-2)}{r} \left( \frac{\Gamma^{\prime}(\frac{r-1}{r})}{\Gamma(\frac{r-1}{r})} - \frac{\Gamma^{\prime}(n\frac{r-1}{r})}{\Gamma(n\frac{r-1}{r})}\right)\bigg)}{(r-1)^{n(n-1)}}
$$

\vskip 4mm

The next propositions describe more  properties of $\Omega_K$. Some were already stated in the introduction.
\vskip 3mm
\noindent
\begin{proposition} \label{prop}
\noindent Let $K$ be a  convex body with centroid at the origin.
\vskip 2mm
(i) \  For all $p >0$,
$$ 
\Omega_{K} \ls \left(\frac{as_{p}(K^{\circ})}{n|K^{\circ}|}\right)^{\frac{n(n+p)}{p}}.
$$
If $K$ is in addition in $C^2_+$, then equality holds if and only if $K$ is an ellipsoid. 
\vskip 2mm
(ii) \  For all $p \geq 0$ 
$$
\Omega_{K} \ls \left(\frac{as_{p}(K)}{n |K^\circ|}\right)^{n+p}. 
$$
If $K$ is in addition in $C^2_+$, then equality holds if and only if $K$ is an ellipsoid. 
\vskip 2mm
(iii) \   $ \Omega_{K}  \ls  \left( \frac{|K|}{|K^{\circ}|}\right)^{n} $. 
If $K$ is in addition in $C^2_+$, then equality holds if and only if $K$ is an ellipsoid. 
\vskip 2mm
(iv) \    $\Omega_{K} \Omega_{K^{\circ}} \ls 1 $.
If $K$ is in addition in $C^2_+$, then equality holds if and only if $K$ is an ellipsoid. 
\end{proposition}
\vskip 2mm
\noindent
{\bf Proof.} 
\par
\noindent
(i)
The first part follows from  Corollary \ref{cor}, Proposition \ref{propWY} (iii) and  the  Remark (ii) after it.
The second part follows from Corollary \ref{cor}, Proposition \ref{propWY} (iii) and  the  Remark (iv) after it.
\vskip 2mm
(ii)  The first part follows from the definition of $\Omega_K$, Proposition \ref{propWY} (ii) and the  Remark (ii) after it. The second part follows from the definition of $\Omega_K$, Proposition \ref{propWY} (ii) and the  Remark (iv) after it.
\vskip 2mm
\noindent
(iii) By (ii), $\Omega_{K} \ls \left(\frac{as_{0}(K)}{n |K^\circ|}\right)^{n} =  \left( \frac{|K|}{|K^{\circ}|}\right)^{n} $.
\vskip 2mm
\noindent
(iv)  is immediate from (iii).

\vskip 4mm
\noindent
\begin{proposition}\label{prop2}
Let $K$ be a  convex body  $\mathbb R^n$ with centroid at the origin.
\vskip 2mm
(i) \   $ \Omega_{K} =\exp{ \left( \frac{1}{|K^{\circ}|} \int_{\partial K^{\circ}} \langle x, N_{K^{\circ}}(x)\rangle \log{ \frac{ \kappa_{K^{\circ}}(x)}{\langle x, N_{K^{\circ}}(x)\rangle^{n+1}} } d\mu_{K^{\circ}}(x) \right) } $.
\vskip 2mm
\noindent
If $K$ is in addition in $C^2_+$, then
\par
(ii) \  $ \Omega_{K} =\exp{ \left( -\frac{1}{|K^{\circ}|} \int_{\partial K} \frac{\kappa_{K}(x)}{\langle x, N_{K}(x)\rangle^{n}} \log{ \frac{ \kappa_{K}(x)}{\langle x, N_{K}(x)\rangle^{n+1}} } d\mu_{K}(x) \right) } $.
\vskip 2mm
\begin{eqnarray*}
&&\hskip -4mm (iii) \  \frac{1}{|K|} \int_{\partial K} \langle x, N_{K}(x)\rangle \log{ \frac{ \kappa_{K}(x)}{\langle x, N_{K}(x)\rangle^{n+1}} } d\mu_{K}(x)  
\  \leq \  
 n \log\frac{|K^{\circ}|}{|K|} \  \leq \  \\
 &&\hskip 50mm \frac{1}{|K^{\circ}|} \int_{\partial K} \frac{\kappa_{K}(x)}{ \langle x, N_{K}(x)\rangle^{n}}  \log{ \frac{\kappa_{K}(x)}{\langle x, N_{K}(x)\rangle^{n+1}} } d\mu_{K}(x).
\end{eqnarray*}

\end{proposition}

\vskip 3mm
\noindent
{\bf Proof.} 
\par
\noindent
(i)
By Corollary \ref{cor}, 
\begin{eqnarray*}
 \log{\Omega_{K}}&=&\log{\left( \lim_{p\rightarrow 0}  \left(\frac{as_{p}(K^{\circ})}{n|K^{\circ}|}\right)^{\frac{n(n+p)}{p}} \right) }= \log{\left( \lim_{p\rightarrow 0}  \left(\frac{as_{p}(K^{\circ})}{n|K^{\circ}|}\right)^{\frac{n^2}{p}} \right) }\\
&=& \lim_{p\rightarrow 0} \frac{n^2}{p} \log{\frac{as_{p}(K^{\circ})}{n|K^{\circ}|} } = n^{2} \lim_{p\rightarrow 0} \frac{ \frac{d}{dp}\left( as_{p}(K^{\circ}) \right)}{as_{p}(K^{\circ})} \\
&=& n^{2} \lim_{p\rightarrow 0} \frac{n(n+p)^{-2} }{  as_{p}(K^{\circ}) } \int_{\partial K^{\circ}} \frac{\kappa_{K^{\circ}}(x)^{\frac{p}{n+p}}}{\langle x, N_{K^{\circ}}(x) \rangle^{\frac{n(p-1)}{n+p}} } \\
& & \hskip 40mm \times \log{ \frac{\kappa_{K^{\circ}}(x)}{\langle x, N_{K^{\circ}}(x) \rangle^{n+1} } } d\mu_{K^{\circ}}(x)  \\
&=& \frac{1}{|K^{\circ}|} \int_{\partial K^{\circ}} \langle x, N_{K^{\circ}}(x)\rangle \log{ \frac{ \kappa_{K^{\circ}}(x)}{\langle x, N_{K^{\circ}}(x)\rangle^{n+1}} } d\mu_{K^{\circ}}(x).
 \end{eqnarray*}
\vskip 2mm
\noindent
(ii) If $K$ is in $C^2_+$, we have by (\ref{def2}) that 
\begin{eqnarray*}
\log{\Omega_{K}} & =& \log{\left( \lim_{p\rightarrow \infty}  \left(\frac{as_{p}(K)}{as_{\infty}(K)}\right)^{n+p} \right) }= \lim_{p \rightarrow \infty} \frac{ \log{ \left( \frac{as_{p}(K)}{as_{\infty}(K)}\right) } }{ (n+p)^{-1}}  \\
&= &
- \lim_{p \rightarrow \infty} \frac{ (n+p)^{2} \frac{d}{dp}\left(as_{p}(K)\right) }{as_{p}(K)} \\
&= &
 - \lim_{p \rightarrow \infty} \frac{ (n+p)^{2}}{as_{p}(K)} \int_{\partial K} \frac{d}{dp} \bigg( \exp \bigg( \log{(\kappa_{K}(x))} \frac{p}{n+p} \\
& &\hskip 35mm  - \log{(\langle x, N_{K}(x) \rangle )} \frac{n(p-1)}{n+p} \bigg) \bigg) d\mu_{K}(x)  \\
 &=&
 -\lim_{p \rightarrow \infty} \frac{ (n+p)^{2}}{as_{p}(K)} \int_{\partial K} \frac{\kappa_{K}(x)^{\frac{p}{n+p}}}{\langle x, N_{K}(x) \rangle^{\frac{n(p-1)}{n+p}} } \bigg( \frac{n}{(n+p)^{2}}\log{( \kappa_{K}(x))} \\
 & &\hskip 40mm - \frac{n(n+1)}{(n+p)^{2}} \log{( \langle x, N_{K}(x) \rangle )} \bigg) d\mu_{K}(x)  \\
& = &
 - \lim_{p \rightarrow \infty} \frac{n}{as_{p}(K)} \int_{\partial K} \frac{\kappa_{K}(x)^{\frac{p}{n+p}}}{\langle x, N_{K}(x) \rangle^{\frac{n(p-1)}{n+p}} } \log{ \frac{\kappa_{K}(x)}{\langle x, N_{K}(x) \rangle^{n+1} } } d\mu_{K}(x)  \\ 
 &= &
 - \frac{n}{as_{\infty}(K)} \int_{\partial K} \frac{\kappa_{K}(x)}{\langle x, N_{K}(x) \rangle^{n}}  \log{ \frac{\kappa_{K}(x)}{\langle x, N_{K}(x) \rangle^{n+1} } } d\mu_{K}(x).
 \end{eqnarray*}

\vskip 2mm
\noindent
(iii)
Combine Proposition \ref{prop} (iii)  with (i) and (ii).

\vskip 4mm
\noindent 
Let $(X, \mu)$ be a measure space  and let  $dP=pd\mu$ and  $dQ=qd\mu$ be probability measures on $X$ that are  absolutely continuous with respect to the measure $\mu$. 
The {\em Kullback-Leibler divergence} or {\em relative entropy} from $P$ to $Q$ is defined as \cite{CT}
\begin{equation}\label{relent}
 D_{KL}(P\|Q)= \int_{X} p\log{\frac{p}{q}} d\mu.
\end{equation}
\vskip 3mm
\noindent 
{\em The information inequality} \cite{CT} holds for the Kullback-Leibler divergence:
Let $P$ and $Q$ be  as above. Then
\begin{equation}\label{Gibbs}
  D_{KL}(P\|Q) \ge 0,
\end{equation}
with equality if and only if $P=Q$.

\vskip 3mm
\noindent 
The invariant $\Omega_K$ is related  to  relative entropies on $K$ and a corresponding information inequality holds, which is exactly the inequality of Proposition \ref{prop} (iii).
\vskip 3mm
\begin{proposition} \label{prop3}
\noindent 
Let $K$ a  convex body in $\mathbb R^n$ that is $C^{2}_{+}$. 
Let 
\begin{equation}\label{PQ}
p(x)= \frac{ \kappa_{K}(x)}{\langle x, N_{K}(x) \rangle^{n} \  n|K^{\circ}|} \, , \   \ q(x)= \frac{\langle x, N_{K}(x) \rangle }{n\  |K|}.
\end{equation}
\noindent 
Then $P=p\  \mu$ and $Q=q \ \mu$ are probability measures on $\partial K$ that are absolutely continuous with respect  to $\mu_{K}$  and
\begin{equation}\label{prop3:eq1} 
D_{KL}(P\|Q) = \log{\left( \frac{|K|}{|K^{\circ}|} \Omega_{K}^{-\frac{1}{n}} \right) }
\end{equation}
and
\vskip 2mm
\noindent
\begin{equation}\label{prop3:eq2} 
D_{KL}(Q\|P) = \log{ \left(  \frac{|K^{\circ}|}{|K|} \Omega_{K^{\circ}}^{-\frac{1}{n}} \right) }.
\end{equation}
\noindent 
Moreover Gibb's inequality implies that
$$ \Omega_{K} \ls \left(\frac{|K|}{|K^{\circ}|} \right)^{n} $$ 
\noindent with equality if and only if $K$ is an ellipsoid.
\end{proposition}

\vskip 3mm
\noindent
{\bf Proof  of Proposition \ref{prop3}.} 
\par
\noindent
As
$$ 
n|K|= \int_{\partial K}  \langle x, N_{K} \rangle d\mu_{K}(x) \, \  \   \mbox{and} \  \ n|K^{\circ}|  = \int_{\partial K}  \frac{ \kappa_{K}(x)}{\langle x, N_{K}(x) \rangle^{n} }  d\mu_{K}(x), 
$$
$\int_{\partial K} p\ d\mu_K= \int_{\partial K} q\  d\mu_K = 1$ and hence $P$ and $Q$ are probability measures that are absolutely continuous with respect  to $\mu_{K}$ on $K$.
\par
(\ref{prop3:eq1}) resp. (\ref{prop3:eq2}) follow from  the definition  of the relative entropy 
(\ref{relent}) and Proposition \ref{prop2} (ii) resp. Proposition \ref{prop2} (i).
\par
By (\ref{Gibbs}), equality holds in the  inequality of the proposition, if and only if for all $x \in \partial K$
$$
\frac{\kappa_{K}(x)}{\langle x, N_{K}(x)\rangle} = \frac{|K|}{|K^{\circ}|} = \mbox{constant}
$$ 
which holds, by the above mentioned  theorem of Petty \cite{Pe}  if and only if $K$ is an ellipsoid.

\vskip 4mm
Let $K$ be a convex body in $\mathbb{R}^n$. Recall that the normalized cone measure $cm_K$
on $\partial K$ is defined as follows:
For every measurable set $A \subseteq \partial K$
\begin{equation}\label{def:conemeas} 
cm_{K}(A)  = \frac{1}{|K|}|\{ta : \ a \in A, t\in [0,1] \}|.
\end{equation}
For more information about cone measures we refer to e.g.,  \cite{BGMN}, \cite{MilGr}, \cite{Na2} and \cite{NaRo}. 
\vskip 3mm
\noindent
The next proposition is well known. It shows that the measures $P$ and $Q$ defined in Proposition \ref{prop3}
 are the cone measures  of $K$  and $K^\circ$. We include the proof for completeness. $N_K:\partial K \rightarrow S^{n-1}$, $x \rightarrow N_K(x)$  is the Gauss map.
\vskip 3mm
\noindent 
\begin{proposition} \label{prop:conemeas}
\noindent 
Let $K$ a  convex body in $\mathbb R^n$ that is $C^{2}_{+}$. Let $P$ and $Q$ be the probability measures on $\partial K$  defined by (\ref{PQ}).  Then
$$
P= N_{K}^{-1}N_{K^{\circ}}cm_{K^{\circ}}\  \   \mbox{and} \   \ Q= cm_{K},
$$
\noindent 
or, equivalently, for every measurable subset $A$ in $ \partial K$
$$ 
P(A)= cm_{ K^{\circ}} \bigg(N_{{K^{\circ}}}^{-1} \big(N_{K} (A)\big)\bigg) \   \ \mbox{and} \   \   Q(A)= cm_{ K}(A).
$$
\end{proposition}

\vskip 2mm
\noindent 
{\bf Proof.}
\noindent 
$$
Q(A)= \frac{1}{n|K|} \int_{A} \langle x, N_{K}(x)\rangle d\mu_K( x) = 
cm_{K}(A).
$$

\noindent 
Also 
$$ P(A)= \int_{A}\frac{ \kappa_{K}(x)}{\langle x, N_{K}(x) \rangle^{n} } \frac{ d\mu_K( x)}{n |K^{\circ}|} =  
\frac{1}{n|K^{\circ}|} \int_{N_{K}(A)} \frac{1}{h_{K}^{n}(u)} d \sigma(u).
$$

\noindent 
Let $B \subseteq \partial K^{\circ}$. Then 
$$cm_{K^{\circ}} (B) = \frac{1}{|K^{\circ}|} \big| \{ x \in \mathbb R^n: \|x\|_{K^{\circ}} \ls 1, \ \frac{x}{\|x\|_{2}} \in N_{K^{\circ}}(B) \}\big| \ .
$$
\noindent 
Let $\Delta =  \{ x\in \mathbb R^n: \|x\|_{K^{\circ}} \ls 1, \ \frac{x}{\|x\|_{2}} \in N_{K^{\circ}}(B) \} \ $. 
We have (see \cite{NaRo})
\begin{eqnarray*}
cm_{ K^{\circ}} (B) &=& \frac{|\Delta|}{|K^{\circ}|} = 
 \frac{1}{|K^{\circ}|} \int_{0}^{\infty} \int_{S^{n-1}} r^{n-1} 1_{\Delta}(r \theta) dr d\sigma(\theta) \\
 &= & \frac{1}{|K^{\circ}|} \int_{N_{K^{\circ}}(B) } \int_{0}^{\frac{1}{\|\theta\|_{K^{\circ}}}} r^{n-1} dr d\sigma(\theta) \\
 &= &
 \frac{1}{n|K^{\circ}|} \int_{N_{K^{\circ}}(B) } \frac{1}{h_{K}^{n}(\theta)} d \sigma(\theta).
 \end{eqnarray*}
\noindent 
Let $B\in \partial K^{\circ}$ be such that $N_{K^{\circ}}(B)=N_{K}(A)$.  This means that $B=N_{K^{\circ}}^{-1}\big(N_{K}(A)\big)$. Then  $P(A) = cm_{K^{\circ}}\bigg(N_{K^{\circ}}^{-1}\big(N_{K}(A) \big)\bigg)$,  which completes the proof. 
  
\vskip 4mm
\noindent
Therefore, with $P$ and $Q$ defined as in (\ref{PQ}),
\begin{equation}\label{relent:conemeas}
 D_{KL}(P\|Q)= D_{KL}\big(N_{K}N_{K^{\circ}}^{-1}cm_{ K^{\circ}}\| cm_{ K}\big)
 \end{equation}
and we get as a corollary to Proposition \ref{prop2} that the invariant $\Omega_K$ is the exponential of the relative entropy of the cone measures of $K$ and $K^\circ$.
\vskip 3mm
\begin{corollary} \label{cor1}
Let $K$ be a  convex body in $C_{+}^{2}$. Then
$$ \Omega_{K}^{1/n} = \frac{|K^{\circ}|}{|K|} \exp{\bigg(- D_{KL}(N_{K}N_{K^{\circ}}^{-1}cm_{ K^{\circ}} \| cm_{K} ) \bigg)} \ . $$
\end{corollary}

\vskip 4mm
\noindent 
Finally, an  isoperimetric inequality holds for the affine invariant $\Omega_K$:
\vskip 3mm
\noindent
\begin{proposition}
\noindent Let $K$ be a convex body in $C_{+}^{2}$ of volume $1$.  Then
$$ \Omega_{K^{\circ}} \ls \Omega_{(\widetilde{B_{n}^{2}})^{\circ}} $$
\noindent with equality if and only if $K= \widetilde{B_{n}^{2}}$.
\end{proposition}
\vskip 2mm
\noindent

{\bf Proof.}
\par
\noindent
The proof follows from the above information inequality for convex bodies together with the  Blaschke Santal\'o inequality and the fact that $\Omega_{(\widetilde{B_{n}^{2}})^{\circ}} = |B_{n}^{2}|^{2n} $.

\vskip 4mm
In the next section we show that the invariant $\Omega_K$ is related to the $L_p$ centroid bodies.

\vskip 5mm

\section{$Z_p(K)$  for $K$ in $C^2_+$}
 
\medskip

In this section we show how $\Omega_K$ is related to the $L_p$ centroid bodies.
The main theorem of the section is Theorem \ref{theorem1}. We assume there that $K$ is symmetric, 
mainly because  the bodies $Z_{p}(K)$ are symmetric by definition. Also, throughout this section we assume that $K$ is of volume $1$.
\vskip 3mm
\noindent
\begin{theorem} \label{theorem1}
Let $K$ be a symmetric convex body in $\mathbb{R}^n$ of volume $1$ that is in $C^2_+$. Then
\vskip 2mm
\begin{eqnarray*} 
\hskip -48mm (i)
\hskip 2mm  \lim _{p \rightarrow \infty} \frac{p}{\log p} \left(|Z_{p}^\circ(K) |-|K^\circ |\right) = \frac{n(n+1)}{2}\  |K^\circ|.
 \end{eqnarray*}
\vskip 2mm
\begin{eqnarray*} 
&&\hskip -13mm (ii) \ \lim _{p \rightarrow \infty}  p \left( |Z_p^\circ(K)| - |K^\circ|- \frac{n(n+1)}{2p} \log p \ |Z_p^\circ(K)| \right) = \\
&& \lim _{p \rightarrow \infty}  p \left( |Z_p^\circ(K)| - |K^\circ|-  \frac{n(n+1)}{2p} \log p \ |K^\circ| \right) = \\
&&\hskip 14mm - \frac{1}{2} \ \int_{S^{n-1}}  h_K(u)^{-n}\  \log \left(2^{n+1} \pi^{n-1}
h_K(u)^{n+1 } f_K(u) \right) d \sigma(u)= \\
&& \hskip 22mm\frac{1}{2} \ \int_{\partial K} \frac{ \kappa(x)}{\langle x, N(x) \rangle^n} \   \log \left(\frac{\kappa(x)}{2^{n+1} \pi^{n-1}
\langle x, N(x) \rangle^{n+1}} \right) d\mu_K(x)
\end{eqnarray*}
\end{theorem}
\vskip 4mm
\noindent 
Thus Theorem \ref{theorem1} shows that if $K$ is a symmetric convex body in $C_{+}^{2}$ of volume $1$, then 
\begin{eqnarray*} 
&&\hskip -10mm \lim_{p \rightarrow \infty} p \left( |Z_{p}^{\circ}(K)| - |K^{\circ}| - \frac{n(n+1)\log{p}}{2p} |Z_{p}^{\circ}(K)| \right) = \\
&&\hskip -10mm  \lim _{p \rightarrow \infty}  p \left( |Z_p^\circ(K)| - |K^\circ|-  \frac{n(n+1)}{2p} \log p \ |K^\circ| \right) = \\
&&\frac{1}{2} \int_{\partial K} \frac{\kappa_{K}(x)}{\langle x, N(x) \rangle^{n}} \log{\left( 2^{n+1}\pi^{n-1}\frac{\kappa_{K}(x)}{\langle x, N(x) \rangle^{n+1} }\right) } d\mu_{K}(x) =\\
&&\frac{\log{\left( 2^{n+1}\pi^{n-1} \right)}}{2} \int_{\partial K} \frac{\kappa_{K}(x)}{\langle x, N(x) \rangle^{n}}  d\mu_{K}(x) \\
&& + \frac{1}{2}\int_{\partial K} \frac{\kappa_{K}(x)}{\langle x, N(x) \rangle^{n}}  \log{\left(\frac{\kappa_{K}(x)}{\langle x, N(x) \rangle^{n+1} }\right) }  d\mu_{K}(x) = \\
&& \log{\left( 2^{n+1}\pi^{n-1} \right)} \frac{ n|K^{\circ}|}{2} -\frac{|K^{\circ}|}{2} \log{\Omega_{K}} =
   -\frac{|K^{\circ}|}{2} \log{\frac{\Omega_{K}}{2^{n(n+1)}\pi^{n(n-1)}}} 
\end{eqnarray*} 
or 
\begin{eqnarray}\label{equ:Zp}
&&\lim_{p \rightarrow \infty} p \left( \frac {|Z_{p}^{\circ}(K)|}{|K^{\circ}|} - \left(1-\frac{n(n+1)\log{p}}{2p}\right) \right) \nonumber\\
&&=  \lim_{p \rightarrow \infty} p \left( (1-\frac{n(n+1)\log{p}}{2p})\frac{|Z_{p}^{\circ}(K)|}{|K^{\circ}|} -1  \right) =  -\frac{1}{2} \log{\frac{\Omega_{K}}{2^{n(n+1)}\pi^{n(n-1)}} }.
\end{eqnarray}
\vskip 3mm
\noindent 
So we have the following
\vskip 3mm
\begin{corollary}
\noindent 
Let $K$ and $C$ be symmetric convex bodies of volume $1$ in $C^{2}_{+}$. Then
\begin{eqnarray*}
&& \hskip -24mm (i) \ 
 \lim_{p \rightarrow \infty} \frac{2p}{n} \left( \frac{(1-\frac{n(n+1)\log{p}}{2p})|Z_{p}^{\circ}(K)|}{|K^{\circ}|} -1 \right)  =\\
 && \hskip -17mm
 \lim_{p \rightarrow \infty} \frac{2p}{n} \left( \frac {|Z_{p}^{\circ}(K)|}{|K^{\circ}|} - \left(1-\frac{n(n+1)\log{p}}{2p}\right) \right)  =
 -\frac{1}{2} \log{\frac{\Omega^\frac{1}{n}_{K}}{2^{n+1}\pi^{n-1}} } \\
&&= (n+1)\log{\left(\frac{2\pi^{\frac{n-1}{n+1}}}{ |K^\circ|}\right)} +  D_{KL}\left(N_{K}N_{K^{\circ}}^{-1}cm_{ K^{\circ}} \| cm_{ K} \right)
\end{eqnarray*}
\begin{eqnarray*}
\hskip 3mm (ii)
\hskip 1mm
\lim_{p \rightarrow \infty} p \left((1-\frac{n(n+1)\log{p}}{2p}) \frac{|Z_{p}^{\circ}(K)|}{|K^{\circ}|} -1 \right) \ge  \frac{1}{2} \log{\left( 2^{n(n+1)}\pi^{n(n-1)} \frac{|K^{\circ}|} {|K|} \right)}. 
\end{eqnarray*}
The corresponding statement for $\lim_{p \rightarrow \infty} p \left( \frac {|Z_{p}^{\circ}(K)|}{|K^{\circ}|} - \left(1-\frac{n(n+1)\log{p}}{2p}\right) \right) $ also holds.
\begin{eqnarray*}
\hskip -9mm (iii)
\hskip 1mm
 \lim_{p \rightarrow \infty} p \left(1-\frac{n(n+1)\log{p}}{2p}\right)\left( \frac{|Z_{p}^{\circ}(K)|}{|K^{\circ}|} -\frac{|Z_{p}^{\circ}(C)|}{|C^{\circ}|}\right) = \frac{1}{2n}\log{\frac{\Omega_{C}}{\Omega_{K}} }.
\end{eqnarray*}

\end{corollary}
\vskip 2mm
\noindent
{\bf Proof.}
\par
\noindent
(i) follows from (\ref{equ:Zp}) and Corollary \ref{cor1}, (ii) follows from Proposition \ref{prop} (iii) and  (iii) follows from (\ref{equ:Zp}).
\vskip 5mm
\noindent
The remainder of the section is devoted to the proof of Theorem \ref{theorem1}.
We need several lemmas and notations. 
\par
Let $x,y>0$. Let $\Gamma(x) = \int_{0}^{\infty} \lambda^{x-1} e^{-\lambda} d\lambda $ be the Gamma function and let 
$ B(x,y) = \int_{0}^{1} \lambda^{x-1} (1-\lambda)^{y-1} d \lambda = \frac{\Gamma(x) \Gamma(y)}{\Gamma(x+y)} $ be the Beta function.
\par
We write $f(p)=g(p)\pm o(p)$, if there exists a function $h(p)$ such that $f(p)=g(p)+ h(p)$ and $\lim_{p\rightarrow \infty} ph(p) =0$, i.e. $h(p)$ has terms of order $\frac{1}{p^2}$ and higher.  Similarly,  $f(p)=g(p)\pm o(p^2)$, if   there exists a function $h(p)$ such that $f(p)=g(p)+ h(p)$ and $\lim_{p \rightarrow \infty} p^2h(p) =0$, i.e. $h(p)$ has terms of order $\frac{1}{p^3}$ and higher. We write$f(p)=g(p)\pm O(p)$, if there exists a function $h(p)$ such that $f(p)=g(p)+ h(p)$ and $\lim_{p\rightarrow \infty} h(p) =0$
\vskip 3mm
\noindent
\begin{lemma} \label{lemma1}
\vskip 2mm
\noindent
 Let $p >0$. Then
\begin{eqnarray*}
(i) \ 
 \left(B\left(p+1, \frac{n+1}{2}\right)\right)^\frac{n}{p} &=& 
 1- \frac{n(n+1)}{2p} \log p + \frac{n}{p} \log \left(\Gamma(\frac{n+1}{2})\right) +\\
 & & \hskip -5mm \frac{n^2(n+1)^2}{8p^2} (\log p )^2 - 
\frac{n^2(n+1)}{2p^2}   \log \left(\Gamma(\frac{n+1}{2})\right) \log p +\\
 && \hskip -5mm \frac{n}{2p^2}\left[n\left(\log \left(\Gamma(\frac{n+1}{2})\right)\right)^2- \frac{n+1}{4}\left(n(n+1)+2(n+3)\right) \right]\\
 &&\hskip -5mm \pm o(p^2).
\end{eqnarray*} 
\vskip 2mm
(ii)  Let $0 \leq a \leq 1$. Then
\begin{eqnarray*}
&&\left( \int_{0}^{1} u^p (1-u)^\frac{n-1}{2} \left( 1-a\left(1-u \right)\right)^\frac{n-1}{2} du 
\right)^\frac{n}{p}= 
 1- \frac{n(n+1)}{2p} \log p + \\
 && \frac{n}{p} \log \left(\Gamma(\frac{n+1}{2})\right) +
\frac{n^2(n+1)^2}{8p^2} (\log p )^2
- \frac{n^2(n+1)}{2p^2}   \log  \left(\Gamma(\frac{n+1}{2})\right) \log p +\\
&&\frac{n}{2p^2}\left[n\left(\log \left(\Gamma(\frac{n+1}{2})\right)\right)^2- \frac{(n+1)\left(n^2+3n+6\right)}{4}
-(n+1)  {\frac{n-1}{2} \choose 1} \ a \right]  \pm o(p^2).
\end{eqnarray*} 
\end{lemma}
\vskip 3mm
\noindent
The proof of Lemma \ref{lemma1} is in the Appendix.

\vskip 4mm
Let $f: \mathbb R_{+} \rightarrow \mathbb R_{+} $ be a $C^{2}$ $\log$-concave function with $\int_{\mathbb R_{+}} f(t) dt < \infty$ and let $p\ge 1$. 
Let $g_p(t) = t^p f(t)$ and let $t_p =t_p(f)$ the unique point such that $g^{'}(t_p)=0$. We make use of the following Lemma due to B. Klartag \cite{Kl1}  (Lemma 4.3 and Lemma 4.5).

\begin{lemma} \label{klar}
\noindent 
Let $f$ be as above.  For every $\varepsilon \in (0,1)$,
$$ \int_{0}^{\infty} t^p f(t) dt \leq \left( 1+ Ce^{-cp\varepsilon^{2}} \right) \int_{t_p(1-\varepsilon)}^{t_p(1+\varepsilon)} t^p f(t) dt $$
where $C>0$ and  $c>0$ are universal constants.
\end{lemma}
We think that the next lemma is well known. We give a proof for completeness.
\vskip 3mm
\noindent 
\begin{lemma} \label{limit=h}
Let $u \in S^{n-1}$. Let $f$ and $t_p$ be as above and $f$ also such that it is  decreasing and a probability density on $[0,h(u)]$. Then
$$
\lim_{p\rightarrow \infty} t_p = h(u). 
$$
\end{lemma}
\vskip 2mm
\noindent{\bf Proof.}
\newline
We only have to show that $\lim_{p\rightarrow \infty} t_p \geq h(u)$. By H\"older,  $\left(\int_{0}^{h(u)} t^p f(t) dt \right)^\frac{1}{p}\rightarrow h(u)$. Thus, for $\varepsilon >0$ given, there exists $p_{\varepsilon}$ such that for all $p \geq p_{\varepsilon}$,
\begin{equation*}
\int_{0}^{h(u)} t^p f(t) dt \geq \big(h(u) -\varepsilon \big)^p
\end{equation*}
By Lemma \ref{klar}, for all $0 < \delta <1 $,  $\int_{0}^{\infty} t^p f(t) dt \leq \left( 1+ Ce^{-cp\delta^{2}} \right) \int_{t_p(1-\delta)}^{t_p(1+\delta)} t^p f(t) dt $. We choose $\delta = \frac{1}{p^\frac{1}{4}}$ with $p>p_{\varepsilon}$ and get, using the monotonicity behavior of $t^pf$ on the respective intervals, that
\begin{eqnarray*}
\big(h(u) -\varepsilon \big)^p & \leq & \left( 1+ Ce^{-cp\sqrt{p}} \right)\left[ \int_{t_p(1-\delta)}^{t_p} t^p f(t) dt + 
\int_{t_p}^{t_p(1+\delta)} t^p f(t) dt \right]\\
& \leq & \left( 1+ Ce^{-cp\sqrt{p}} \right)\  p^\frac{1}{4} t_p f(t_p) \ t_p^p.
\end{eqnarray*}
As $f$ is decreasing, $f(t_p) \leq f(0)$. Moreover $t_p \leq h(u)$. Thus, for $p \geq p_{\varepsilon}$ large enough,
$\left(p^\frac{1}{4} t_p f(t_p)\right)^\frac{1}{p} \leq 1+\varepsilon$ and hence $h(u) - \varepsilon < (1+\varepsilon) \ t_p$

\vskip 3mm
\noindent 
{\bf Remark}
\par
\noindent 
We will apply Lemma \ref{klar} to the function $f(t)=|K\cap (u^{\perp} +t u)|$, $u \in S^{n-1}$. We show below that $f$ is $C^2$.  Thus $t_p$ is well defined and    Lemma \ref{klar} holds. 
Also,  $t_p$ is an increasing function of $p$ and by the above Lemma, \ref{limit=h}, $\lim_{p\rightarrow \infty} t_p = h_K(u)$. 

\vskip 3mm
We also think that the following lemma is well known but we could not find a proof in the literature.
Therefore we include a proof. 
\vskip 2mm
\begin{lemma} \label{derivative}
\noindent 
Let $K$ be a convex body in in $C^2_+$. Let $u \in S^{n-1}$ and let $H_t$ be the hyperplane
orthogonal to $u$ at distance $t$ from the origin. Let $f(t) = |K \cap H_t|$. 
Then $f$ is $C^2$. In fact, 
$$
f^\prime (t) = -  \int_{\partial K \cap H_t}  \frac{\langle u, N_K(x) \rangle} 
{\big(1- \langle u, N_K(x) \rangle^2 \big)^\frac{1}{2}}\  d \mu_{\partial K \cap H_t}(x)
$$
and
$$
f^{\prime \prime} (t) = 
$$
$$
-  \int_{\partial K \cap H_t} \bigg[ \frac{\kappa(x_t)^\frac{1}{n-1}}{ \left(1- \langle N_K(x_t), u \rangle ^2\right)^\frac{3}{2}} - \frac{(n-2) \   \langle N_K(x_t), u \rangle^2 }{\langle N_{K \cap H_t}(x_t), x_t \rangle \   \left(1- \langle N_K(x_t), u \rangle ^2\right)}\bigg]  d \mu_{\partial K \cap H_t}(x_t).
$$

\end{lemma}

\vskip 3mm
\noindent
{\bf Proof.}
\newline
We assume that $\mbox{int}(K) \cap H_t \neq \emptyset$.  To show that $f \in C^2$, we compute
the derivates of  $f$. We first show that 
\begin{eqnarray*}
f^\prime (t) = -  \int_{\partial K \cap H_t}  \frac{\langle u, N_K(x) \rangle} 
{\big(1- \langle u, N_K(x) \rangle^2 \big)^\frac{1}{2}}\  d \mu_{\partial K \cap H_t}(x).
\end{eqnarray*}
Indeed, for $x \in \partial K \cap H_t$, let $\alpha (x) $ be the (smaller) angle formed by $N_K(x)$ and $u$.
Then $\mbox{cos}\  \alpha (x) = \langle u, N_K(x) \rangle$ and 
\begin{eqnarray*}
f^\prime (t)  =  \lim_{\varepsilon \rightarrow 0} \frac{1}{\varepsilon} \  \bigg( |K \cap H_{t+ \varepsilon}| - |K \cap H_t| \bigg) 
 =   -  \lim_{\varepsilon \rightarrow 0} \frac{1}{\varepsilon} \   \bigg( \int_{\partial K \cap H_t} \varepsilon \ \mbox{cot} \ \alpha (x)  \ d \mu_{\partial K \cap H_t}(x) \bigg) \\
 =  -  \int_{\partial K \cap H_t}  \frac{\langle u, N_K(x) \rangle} 
{\big(1- \langle u, N_K(x) \rangle^2 \big)^\frac{1}{2}}\  d \mu_{\partial K \cap H_t}(x).
\end{eqnarray*}
We show next that 
$$
f^{\prime \prime} (t) = 
$$
$$
-  \int_{\partial K \cap H_t} \bigg[ \frac{\kappa(x_t)^\frac{1}{n-1}}{ \left(1- \langle N_K(x_t), u \rangle ^2\right)^\frac{3}{2}} - \frac{(n-2) \   \langle N_K(x_t), u \rangle^2 }{\langle N_{K \cap H_t}(x_t), x_t \rangle \   \left(1- \langle N_K(x_t), u \rangle ^2\right)}\bigg]  d \mu_{\partial K \cap H_t}(x_t).
$$
\par
\noindent
By definition
\begin{eqnarray*} 
f^{\prime \prime} (t)  = - \lim_{\varepsilon \rightarrow 0} \  \frac{1}{\varepsilon} \  \bigg( 
 \int_{\partial K \cap H_{t+\varepsilon}}  \frac{\langle u, N_K(y_{t+\varepsilon}) \rangle} 
{\big(1- \langle u, N_K(y_{t+\varepsilon}) \rangle^2 \big)^\frac{1}{2}} \  d \mu_{\partial K \cap H_{t + \varepsilon}}(y_{t+\varepsilon}) \\
\hskip 40mm - \int_{\partial K \cap H_t}  \frac{\langle u, N_K(x_t) \rangle} 
{\big(1- \langle u, N_K(x_t) \rangle^2 \big)^\frac{1}{2}}\  d \mu_{\partial K \cap H_t}(x_t) 
\bigg)
\end{eqnarray*}
We  project  $ K \cap H_{t+\varepsilon}$ onto $K \cap H_{t} $ and we want to integrate both expressions over 
$\partial K \cap H_{t} $. 
To do so, we fix -  after the projection -  an interior point $x_0$ in $ K \cap H_{t + \varepsilon} $.
For $x_t \in \partial K \cap H_t$, let $[x_0, x_t]$ be the line segment from $x_0$ to $x_t$ and  let $x_{t + \varepsilon}= \partial K \cap H_{t + \varepsilon} \cap [x_0, x_t]$. Now observe that
$$
d \mu_{\partial K \cap H_{t + \varepsilon}} =  \frac{1}{\langle N_{K \cap H_t}(x_t), N_{K \cap H_{t + \varepsilon}}(x_{t + \varepsilon} ) \rangle} \  \left(\frac{\|x_{t + \varepsilon} \|}{\|x_{t }\|}\right)^{n-2} \        d \mu_{\partial K \cap H_t},
$$
where $N_{K \cap H_t}(x_t)$ is the outer normal in $x_t$ to the boundary of the $n-1$ dimensional convex body $K \cap H_t$ and similarly,  $N_{K \cap H_{t + \varepsilon} }(x_{t + \varepsilon})$ is the outer normal in $x_{t + \varepsilon}$ to the boundary of the $n-1$ dimensional convex body $K \cap H_{t + \varepsilon}$.
\par
Notice further that 
$$
\|x_t\| - \|x_{t + \varepsilon} \| = \frac{\varepsilon \   \langle N_K(x_t), u \rangle  \   \|x_t\| }{\langle N_{K \cap H_t}(x_t), x_t \rangle \   \left(1- \langle N_K(x_t), u \rangle ^2\right)^\frac{1}{2}} + \mbox{ higher order terms in} \  \varepsilon.
$$
Therefore
\begin{eqnarray*}
\left(\frac{\|x_{t + \varepsilon} \|}{\|x_{t }\|}\right)^{n-2} &=& \left(1- \frac{\varepsilon \   \langle N_K(x_t), u \rangle }{\langle N_{K \cap H_t}(x_t), x_t \rangle \   \left(1- \langle N_K(x_t), u \rangle ^2\right)^\frac{1}{2}}  \right) ^{n-2} \\
&=& 1-  \frac{(n-2)\  \varepsilon \   \langle N_K(x_t), u \rangle }{\langle N_{K \cap H_t}(x_t), x_t \rangle \   \left(1- \langle N_K(x_t), u \rangle ^2\right)^\frac{1}{2}} \\
&+& \mbox{ higher order terms in} \  \varepsilon
\end{eqnarray*}
Thus
\begin{eqnarray*} 
f^{\prime \prime} (t) & = &-\  \lim_{\varepsilon \rightarrow 0} \  \frac{1}{\varepsilon} \   
\int_{\partial K \cap H_t} \bigg[\frac{\langle u, N_K(y_{t+\varepsilon}) \rangle} 
{\langle N_{K \cap H_t}(x_t), N_{K \cap H_{t + \varepsilon}}(x_{t + \varepsilon} ) \rangle\  \big(1- \langle u, N_K(y_{t+\varepsilon}) \rangle^2 \big)^\frac{1}{2}} \\
&\times& \bigg( 1-  \frac{(n-2)\  \varepsilon \   \langle N_K(x_t), u \rangle }{\langle N_{K \cap H_t}(x_t), x_t \rangle \   \left(1- \langle N_K(x_t), u \rangle ^2\right)^\frac{1}{2}}
+ \mbox{ higher order terms in} \  \varepsilon \bigg) \\
&- &\  \frac{\langle u, N_K(x_t) \rangle} 
{\big(1- \langle u, N_K(x_t) \rangle^2 \big)^\frac{1}{2}} \bigg]\  d \mu_{\partial K \cap H_t}(x_t) \\
&=&- \int_{\partial K \cap H_t} \lim_{\varepsilon \rightarrow 0} \  \frac{1}{\varepsilon} \ 
\bigg[\frac{\langle u, N_K(y_{t+\varepsilon}) \rangle} 
{\langle N_{K \cap H_t}(x_t), N_{K \cap H_{t + \varepsilon}}(x_{t + \varepsilon} ) \rangle\  \big(1- \langle u, N_K(y_{t+\varepsilon}) \rangle^2 \big)^\frac{1}{2}} \\
&\times& \bigg( 1-  \frac{(n-2)\  \varepsilon \   \langle N_K(x_t), u \rangle }{\langle N_{K \cap H_t}(x_t), x_t \rangle \   \left(1- \langle N_K(x_t), u \rangle ^2\right)^\frac{1}{2}}
+ \mbox{ higher order terms in} \  \varepsilon \bigg) \\
&- &\  \frac{\langle u, N_K(x_t) \rangle} 
{\big(1- \langle u, N_K(x_t) \rangle^2 \big)^\frac{1}{2}} \bigg]\  d \mu_{\partial K \cap H_t}(x_t).
\end{eqnarray*}
We can interchange integration and limit using Lebegue's theorem as the functions under the integral are uniformly (in $t$)  bounded by a constant. 
\par
Denote $g_x(t)= \frac{\langle N_K(x_t), u \rangle}{\big(1- \langle u, N_K(x_t) \rangle^2 \big)^\frac{1}{2}}$. Then the expression under the integral becomes
\begin{eqnarray*} 
\lim_{\varepsilon \rightarrow 0} \  \frac{1}{\varepsilon}   &&\hskip -4mm \bigg[ \frac{g_y(t+\varepsilon)}{\langle N_{K \cap H_t}(x_t), N_{K \cap H_{t + \varepsilon}}(x_{t + \varepsilon})\rangle } \bigg( 1-  \frac{(n-2)\  \varepsilon \   \langle N_K(x_t), u \rangle }{\langle N_{K \cap H_t}(x_t), x_t \rangle \   \left(1- \langle N_K(x_t), u \rangle ^2\right)^\frac{1}{2}} \\
&& + \mbox{ higher order terms in} \  \varepsilon \bigg)  -g_x(t) \bigg] \\
&&=  \lim_{\varepsilon \rightarrow 0} \  \frac{1}{\varepsilon} \  \bigg[ g_y(t+\varepsilon) - g_x(t)\bigg] - \frac{(n-2) \   \langle N_K(x_t), u \rangle^2 }{\langle N_{K \cap H_t}(x_t), x_t \rangle \   \left(1- \langle N_K(x_t), u \rangle ^2\right)}.
\end{eqnarray*} 
Here we have also  used that,  as $\varepsilon \rightarrow 0$, $x_{t + \varepsilon} \rightarrow x_t$,  
$N_{K \cap H_{t + \varepsilon}}(x_{t + \varepsilon}) \rightarrow N_{K \cap H_{t}}(x_t)$
and $g_y(t +\varepsilon) \rightarrow g_x(t)$.
\par
To compute $\lim_{\varepsilon \rightarrow 0} \  \frac{1}{\varepsilon} \  \bigg[ g_y(t+\varepsilon) - g_x(t)\bigg] $, we approximate the boundary of $\partial K$ in   $x_t$ by an ellipsoid. This can be done as $\partial K$ is $C^2_+$ by assumption (see Lemma \ref{Dupin} below). To simplify the computations, we assume that the approximating ellipsoid is a Euclidean ball. The case of the ellipsoid is treated similarly,  the computations are just slightly more involved.
As the expression under the integral depends only on the angles between the vectors involved, we can put the origin so that  the approximating Euclidean ball is centered at $0$.  Let $r= \kappa(x_t)^\frac{-1}{n-1}$ be its radius. Then
\begin{eqnarray*}
\lim_{\varepsilon \rightarrow 0} \  \frac{1}{\varepsilon} \  \bigg[ g_y(t+\varepsilon) - g_x(t)\bigg] = \frac{1}{r \  \left(1- \langle N_K(x_t), u \rangle ^2\right)^\frac{3}{2}} = \frac{\kappa(x_t)^\frac{1}{n-1}}{ \left(1- \langle N_K(x_t), u \rangle ^2\right)^\frac{3}{2}}.
\end{eqnarray*}
Alltogether
\begin{eqnarray*}
f^{\prime \prime} (t) = 
-  \int_{\partial K \cap H_t} \bigg[ \frac{\kappa(x_t)^\frac{1}{n-1}}{ \left(1- \langle N_K(x_t), u \rangle ^2\right)^\frac{3}{2}} - \frac{(n-2)    \langle N_K(x_t), u \rangle^2 }{\langle N_{K \cap H_t}(x_t), x_t \rangle    \left(1- \langle N_K(x_t), u \rangle ^2\right)}\bigg]  d \mu_{\partial K \cap H_t}(x_t).
\end{eqnarray*}

\vskip 4mm
\noindent
\begin{lemma} \label{interchange}
Let $K$ be a symmetric convex body of volume $1$ in $C^2_+$.
\vskip 2mm
(i) The functions 
$$
\frac{p}{\mbox{log}(p) }  \frac{1}{h_{Z_p(K)}(u)^n} \left(1-\frac{h_{Z_p(K)}(u)^n}{h_{K}(u)^n}\right)
$$
are uniformly (in $p$) bounded by a function that is integrable on $S^{n-1}$.
\vskip 2mm
(ii) The functions 
$$
\frac{p}{h_{Z_p(K)}(u)^n}   \left(1-\frac{h_{Z_p(K)}(u)^n}{h_{K}(u)^n} - \frac{n(n+1)}{2}\  
\frac{\mbox{log}(p)}{p}\  \frac{h_{Z_p(K)}(u)^n}{h_{K}(u)^n}\right)
$$
are uniformly (in $p$) bounded by a function that is integrable on $S^{n-1}$.
\end{lemma}

\vskip 2mm
\noindent
{\bf Proof.}
\newline
(i) Let $u \in S^{n-1}$. Let $x \in \partial K$ be such that $N_K(x)=u$. As $K$ is in $C^2_+$, by the
Blaschke rolling theorem (see \cite{Sch}), there exists a ball with radius $r_0$ that rolls freely in $K$:
for all $x \in \partial K$,
$B^n_2(x-r _0N(x), r_0) \subset K$.
As $K$ is symmetric, 
\begin{eqnarray*}
h_{Z_p}(u)^n &=& \left(2 \int_0^{h_K(u)}  t^p |\{y \in K: \langle u,y \rangle =t\}| dt \right)^\frac{n}{p} \\
&\geq &  \left( 2\ \int_{h_K(u)-r}^{h_K(u)}  t^p | 
\{y \in B^n_2\big(x- r_0\ u, r_0\big)  : \langle u,y \rangle =t\}|\ dt \right)^\frac{n}{p}\\
&=&
2^\frac{n}{p}\ |B^{n-1}_2|^\frac{n}{p}  \left( \int_{h_K(u)-r_0}^{h_K(u)}  t^p
\bigg( 2 r_0 \big(h_K(u)-t\big) \left[ 1-\frac{h_K(u)-t}{2r_0}\right]\bigg)^\frac{n-1}{2} dt \right)^\frac{n}{p}
\end{eqnarray*}
\par
\noindent
The equality holds as the $(n-1)$-dimensional Euclidean ball 
$$
B^n_2 (x-r_0\ u, r) \cap \{y \in \mathbb{R}^n: \langle u, y \rangle =t\}
$$
has radius $\bigg( 2 r_0\big(h_K(u)-t\big) \left[ 1-\frac{h_K(u)-t}{2  r_0}\right]\bigg)^\frac{1}{2}$.
Now  -  where, to abbreviate, we write $h_K$,  $h_{Z_p(K)}$, instead of $h_K(u)$, $h_{Z_p(K)}(u)$ - 
and  where we use that $ \frac{1}{2} \leq 1-\frac{h_K(u)-t}{2r_0}$,
\begin{eqnarray} \label{est: below}
h_{Z_p}(u)^n &\geq & 2^\frac{n}{p}\ |B^{n-1}_2|^\frac{n}{p}  \big(r_0\ h_K\big)^\frac{n(n-1)}{2p} 
\left( \int_{h_K-r_0}^{h_K}  t^p \left(1-\frac{t}{h_K}\right)^\frac{n-1}{2} dt 
\right)^\frac{n}{p}\nonumber  \\ 
&= &h_K^n \  \left(2\ |B^{n-1}_2| \ h_K^\frac{n+1}{2}\ r_0^\frac{n-1}{2} \right)^\frac{n }{p} \   \left( \int_{1-\frac{r_0}{h_K}}^{1} w^p (1-w)^\frac{n-1}{2} dw  \right)^\frac{n}{p}.
\end{eqnarray}
As $K$ is symmetric, $ r_0 \leq h_K(u)$. If $r_0 = h_K(u)$, then
\begin{eqnarray*}
\frac{h_{Z_p(K)}^n}{h_{K}^n} 
\geq  \left( 2  \ r_0^\frac{n-1}{2} \  h_K^\frac{n+1}{2}\  |B^{n-1}_2|\right)^\frac{n}{p} 
\left( \int_{0}^{1} w^p (1-w)^\frac{n-1}{2}dw  \right)^\frac{n}{p}.
\end{eqnarray*}
If $r_0 < h_K(u)$, we apply Lemma \ref{klar} to the function $f(w)=(1-w)^\frac{n-1}{2}$. We choose 
$\varepsilon$ so small and $p_0$ so large that $\varepsilon +(1+\varepsilon)\frac{n-1}{2p_0} \leq \frac{r_0}{h_K}$. Then Lemma \ref{klar} holds and we get for all $p \geq p_0$
\begin{eqnarray*}
\frac{h_{Z_p(K)}^n}{h_{K}^n} 
&\geq&  \left( 2  \ r_0^\frac{n-1}{2} \  h_K^\frac{n+1}{2}\  |B^{n-1}_2|\right)^\frac{n}{p} 
\left( \int_{1-\frac{r_0}{h_K}}^{1} w^p (1-w)^\frac{n-1}{2} dw 
\right)^\frac{n}{p}\\
&\geq &
\left( \frac{2\   r_0^\frac{n-1}{2} \   h_K^\frac{n+1}{2}\  |B^{n-1}_2|}{1+C\ e^{-cp\varepsilon^2}}\right)^\frac{n}{p} 
\left( \int_{0}^{1} w^p (1-w)^\frac{n-1}{2} dw 
\right)^\frac{n}{p}\\
&= &\left(   \frac{2\   r_0^\frac{n-1}{2} \   h_K^\frac{n+1}{2}\  |B^{n-1}_2|}{1+C\ e^{-cp\varepsilon^2}}\right)^\frac{n}{p} 
\left( B(p+1, \frac{n+1}{2})\right)^\frac{n}{p}.
\end{eqnarray*}
As
$$
 \left(  2  \ r_0^\frac{n-1}{2} \  h_K^\frac{n+1}{2}\  |B^{n-1}_2|\right)^\frac{n}{p} =
1+ \frac{n}{p} \log\left[ 2  \ r_0^\frac{n-1}{2} \  h_K^\frac{n+1}{2}\  |B^{n-1}_2|\right] \pm o(p),
$$
respectively
$$
\left( \frac{2\   r_0^\frac{n-1}{2} \   h_K^\frac{n+1}{2}\  |B^{n-1}_2|}{1+C\ e^{-cp\varepsilon^2}}\right)^\frac{n}{p} = 1+ \frac{n}{p} \log\left[ \frac{2\   r_0^\frac{n-1}{2} \   h_K^\frac{n+1}{2}\  |B^{n-1}_2|}{1+C\ e^{-cp\varepsilon^2}}\right] \pm o(p)
$$
we get,  together with Lemma \ref{lemma1} (i) 
\begin{eqnarray}\label{eq:inter1}
 \frac{h_{Z_p(K)}^n}{h_{K}^n}  \nonumber 
 &\geq &1- \frac{n(n+1)}{2p} \log p + \frac{n}{p} \  
\log\left[ 2  \ r_0^\frac{n-1}{2} \  h_K^\frac{n+1}{2}\  |B^{n-1}_2| \ \Gamma\left(\frac{n+1}{2}\right)\right] \pm o(p)\\
& \geq& 1- \frac{n(n+1)}{2p} \log p + \frac{n}{2p} \  
\log \left[ 4\  r_0^{n-1} \pi^{n-1} h_K^{n+1} \right] \pm o(p)
\end{eqnarray}
respectively
\begin{eqnarray}\label{eq:inter2}
\frac{h_{Z_p(K)}^n}{h_{K}^n}  
& \geq& 1- \frac{n(n+1)}{2p} \log p + \frac{n}{p} \  
\log\left[ \frac{2\   r_0^\frac{n-1}{2} \   h_K^\frac{n+1}{2}\  |B^{n-1}_2| \ \Gamma\left(\frac{n+1}{2}\right)}{1+C\ e^{-cp\varepsilon^2}}\right] \pm o(p)\\
& \geq& 1- \frac{n(n+1)}{2p} \log p + \frac{n}{2p} \  
\log\left[ \frac{4 \ r_0^{n-1} \pi^{n-1} h_K^{n+1}}{(1+C\ e^{-cp\varepsilon^2})^2} \right] \pm o(p)
\end{eqnarray}

\noindent
Now notice that there is $\alpha >0$ such that 
$$
B^n_2(0,  \alpha) \subset K \subset B^n_2(0, \frac{1}{\alpha}).
$$
This implies that for all $u \in S^{n-1}$
$\alpha \leq h_K \leq \frac{1}{\alpha}$. Moreover we can choose $\alpha$ so small that we have for all $p \geq p_0 >1$
$$
B^n_2(0,  \alpha) \subset Z_p(K)  \subset K \subset B^n_2(0,  \frac{1}{\alpha}),
$$
which implies that for all $u \in S^{n-1}$, for all $p \geq p_0$, 
\begin{equation}\label{alpha}
\alpha \leq h_{Z_p(K)} \leq \frac{1}{\alpha}.
\end{equation} 
On the one hand, as $Z_p(K)  \subset K$,
$$
\frac{p}{\mbox{log}(p) }  \frac{1}{h_{Z_p(K)}(u)^n} \left(1-\frac{h_{Z_p(K)}(u)^n}{h_{K}(u)^n}\right)
\geq 0
$$
On the other hand, we get by (\ref{eq:inter1}), (\ref{eq:inter2}) and (\ref{alpha}) with a constant $c$
\begin{eqnarray*} 
\frac{p}{\mbox{log}(p) }  \frac{1}{h_{Z_p(K)}(u)^n} \left(1-\frac{h_{Z_p(K)}(u)^n}{h_{K}(u)^n}\right)
&\leq &
\frac{ c n}{ \alpha^n}\left( n+1 - \frac{1}{\log p} 
\log \left(4 r_0^{n-1} \pi^{n-1} h_K^{n+1} \right)\right)\\
&\leq &
\frac{ c n}{ \alpha^n}\left( n+1 + \frac{1}{\log p_0} 
\bigg|\log \left(\frac{4 r_0^{n-1} \pi^{n-1} }{\alpha^{n+1}} \right)\bigg| \right)
\end{eqnarray*}
respectively
\begin{eqnarray*}
\frac{p}{\mbox{log}(p) }  \frac{1}{h_{Z_p(K)}(u)^n} \left(1-\frac{h_{Z_p(K)}(u)^n}{h_{K}(u)^n}\right)
&\leq &
\frac{ c n}{ \alpha^n}\left( n+1 - \frac{1}{\log p} 
\log \left(\frac{4 r_0^{n-1} \pi^{n-1} h_K^{n+1}}{(1+C\ e^{-cp\varepsilon^2})^2} \right)\right)\\
&\leq &
\frac{ c n}{ \alpha^n}\left( n+1 + \frac{1}{\log p_0} 
\bigg|\log \left(\frac{4 r_0^{n-1} \pi^{n-1} }{\alpha^{n+1}} \right) \bigg| \right)
\end{eqnarray*}

The right hand side is a constant and hence integrable.

\vskip 3mm
\noindent
(ii)
As $K$ is in $C^2_+$, there is $R \geq r_0 > 0$ such that for all $x \in \partial K$, $K \subset B^n_2(x-R N(x), R) $. Then we show similarly to (\ref{est: below})  that
\begin{eqnarray*}
h_{Z_p}(u)^n &\leq & h_K^n \  \left(2^\frac{n-1}{2} \ |B^{n-1}_2| \ h_K^\frac{n+1}{2}\ R^\frac{n-1}{2} \right)^\frac{n }{p} \   \left( \int_{0}^{1} w^p (1-w)^\frac{n-1}{2} dw  \right)^\frac{n}{p}.
\end{eqnarray*}
and thus, similar to  (\ref{eq:inter1})
\begin{eqnarray*}
 \frac{h_{Z_p(K)}^n}{h_{K}^n}  
 \leq 1- \frac{n(n+1)}{2p} \log p + \frac{n}{2p} \  
\log \left[ 2^{n+1} \  R^{n-1} \pi^{n-1} h_K^{n-1} \right] \pm o(p)
\end{eqnarray*}
\noindent
Hence, together with  (\ref{eq:inter1}) respectively (\ref{eq:inter2})
\begin{eqnarray*} 
&&\hskip -11mm -\frac{n}{2\ h_{Z_p(K)^n} }\  
\log\left[ 2^{n+1} \  R^{n-1} \pi^{n-1}\    h_K^{n-1} \right]  
\pm O(p) \leq \\
&&\frac{p}{h_{Z_p(K)}(u)^n}  \left(1-\frac{h_{Z_p(K)}(u)^n}{h_{K}(u)^n} - \frac{n(n+1)}{2}\  
\frac{\mbox{log}(p)}{p}\  \frac{h_{Z_p(K)}(u)^n}{h_{K}(u)^n} \right)\\
&& \hskip 45mm \leq 
-\frac{n}{2\ h_{Z_p(K)^n} }\  
\log\left[ 4 \  r_0^{n-1} \pi^{n-1}\    h_K^{n+1} \right]  
\pm O(p).
\end{eqnarray*} 
respectively
\begin{eqnarray*} 
&&\hskip -20mm -\frac{n}{2\ h_{Z_p(K)^n} }\  
\log\left[ 2^{n+1} \  R^{n-1} \pi^{n-1}\    h_K^{n-1} \right]  
\pm O(p) \leq \\
&&\frac{p}{h_{Z_p(K)}(u)^n}  \left(1-\frac{h_{Z_p(K)}(u)^n}{h_{K}(u)^n} - \frac{n(n+1)}{2}\  
\frac{\mbox{log}(p)}{p}\  \frac{h_{Z_p(K)}(u)^n}{h_{K}(u)^n} \right)\\
&& \hskip 45mm \leq 
-\frac{n}{2\ h_{Z_p(K)^n} }\  
\log\left[ \frac{4 \  r_0^{n-1} \pi^{n-1}\    h_K^{n+1} }{(1+C\ e^{-cp\varepsilon^2})^2}\right]  
\pm O(p).
\end{eqnarray*} 

\noindent
Hence, using (\ref{alpha}), we get with an  absolute constant $c$ for all $p \geq p_0$ 
\begin{eqnarray*} 
\bigg|  \   \frac{p}{h_{Z_p(K)}(u)^n}  \left(1-\frac{h_{Z_p(K)}(u)^n}{h_{K}(u)^n} - \frac{n(n+1)}{2}\  
\frac{\mbox{log}(p)}{p}\  \frac{h_{Z_p(K)}(u)^n}{h_{K}(u)^n}\right)  \  \bigg| \\
\leq \frac{c n }{\alpha^n}  
\bigg|  \log\left[ \frac{2^{n+1} \  R^{n-1} \pi^{n-1}}{\alpha^{n-1} }\right]  \bigg| 
\end{eqnarray*} 

Again, the right hand side is a constant and therefore integrable.

\vskip 4mm
As $K \in C^2_+$, the indicatrix of Dupin at every $x \in \partial K$ is an ellipsoid. 
Since the quantities  considered in the above Theorem \ref{theorem1} are  affine invariant,  we can assume that the indicatrix
is a Euclidean ball.
We have (see \cite{SW2})
\vskip 3mm
\noindent
\begin{lemma} \label{Dupin}
\noindent 
Let $K \subset \mathbb R^n$ be a convex body in $C^2_+$.  
We assume that the indicatrix of Dupin at $x \in \partial K$ is a Euclidean ball.  Let 
$r=r(x)= \kappa(x)^{-\frac{1}{n-1}}$  and put $u=N_{K} (x)$. $B(x-ru,r)$ 
is the Euclidean ball with center at $x-ru$ and radius $r$. Then
for every  $\varepsilon >0$ there exists $\Delta_\varepsilon>0$ such that for all $\Delta\leq \Delta_\varepsilon$,
$$
B\big(x-(1-\varepsilon) r u, (1-\varepsilon)r \big) \cap H (x-\Delta u, u) ^- 
$$
$$
\subset K  \cap  H (x-\Delta u, u) ^- \subset B\big(x-(1+\varepsilon) r u, (1+\varepsilon) r\big) \cap  H (x-\Delta u, u) ^-.
$$
\end{lemma}
$H (x-\Delta u, u) $ is the hyperplane  with normal $u$ through $x-\Delta u$ and $H (x-\Delta u, u)^- $
is the half space determined by this hyperplane  into which $u$ points.
\vskip 4mm
\noindent
{\bf Proof of Theorem \ref{theorem1}}
\vskip 3mm
\noindent
(i) 
$$
|Z_{p}^\circ(K) |-|K^\circ |= \frac{1}{n}\  \int_{S^{n-1}} \left(\frac{1}{h_{Z_p(K)}^{n}(u)} - \frac{1}{h_{K}^{n}(u)}\right) d \sigma(u)
$$
Hence
\begin{eqnarray*}
&& \hskip -5mm 
 \lim _{p \rightarrow \infty} \frac{p}{\log p} \left(|(Z_{p}^\circ(K)) |-|K^\circ |\right)
 =  \frac{1}{n} \lim _{p \rightarrow \infty} \frac{p}{\log p} 
\  \int_{S^{n-1}} \frac{1}{h_{Z_p(K)}^{n}(u)} \left(1 - \frac{h_{Z_p(K)}^{n}(u)}{h_{K}^{n}(u)}\right) d \sigma(u)\\
&&\hskip 40mm = \frac{1}{n}  \int_{S^{n-1}} \lim _{p \rightarrow \infty} \frac{p}{\log p}\  \frac{1}{h_{Z_p(K)}^{n}(u)} \left(1 - \frac{h_{Z_p(K)}^{n}(u)}{h_{K}^{n}(u)}\right) d \sigma(u),
\end{eqnarray*}
where we have used Lemma \ref{interchange} (i) and Lebegue's theorem to interchange integration and limit.
Let $u \in S^{n-1}$. Let $x \in \partial K$ be such that $N_K(x)=u$. As $K$ is in $C^2_+$, $\kappa=\kappa_K(x) >0$ and we can assume that the indicatrix of Dupin at $x$
is a Euclidean ball with radius $r=r(x)= \kappa(x)^\frac{-1}{n-1}$.
\begin{eqnarray*}
&& h_{Z_p(K)}^{n}(u) = \left(\int_K |\langle y,u \rangle|^p dy\right)^\frac{n}{p} 
 = 
 \left(2 \int_0^{h_K(u)}  t^p |\{y \in K: \langle u,y \rangle =t\}| dt \right)^\frac{n}{p} \\
&&\geq  \left(2 \int_{(1-\varepsilon) (h_K(u)- \Delta_\varepsilon)}^{h_K(u)}  t^p |\{y \in K: \langle u,y \rangle =t\}| \ dt \right)^\frac{n}{p} \\
&&\geq   \left( 2\ \int_{(1-\varepsilon) (h_K(u)- \Delta_\varepsilon)}^{h_K(u)}  t^p | 
\{y \in B \big(x-(1-\varepsilon) r\ u, (1-\varepsilon) r \big)  : \langle u, y \rangle =t \}| dt \right)^\frac{n}{p},
\end{eqnarray*}
where we have applied Lemma \ref{Dupin}. In addition, we also choose $\Delta_\varepsilon$ of Lemma \ref{Dupin}  so that $\Delta_\varepsilon \leq \mbox{min}\{\varepsilon, (1-\varepsilon) r\}$.
\vskip 2mm
\noindent
$B(x-(1-\varepsilon) r\ u, (1-\varepsilon) r) \cap \{y \in \mathbb{R}^n: \langle u,y \rangle =t\}$ is a 
$(n-1)$-dimensional Euclidean ball  with radius
$$\bigg( 2 (1-\varepsilon) r \big(h_K(u)-t\big) \left[ 1-\frac{h_K(u)-t}{2 (1-\varepsilon) r}\right]\bigg)^\frac{1}{2}, $$
which, by choice of $\Delta_\varepsilon$ is  bigger or equal than
$$\bigg( 2 (1-\varepsilon) r \big(h_K(u)-t\big) \left[ 1-\frac{\varepsilon\left(h_K(u)+1 -\varepsilon\right)}{2 (1-\varepsilon) r}\right]\bigg)^\frac{1}{2}.$$
Hence
\begin{eqnarray*}
 &&h_{Z_p(K)}^{n}(u) = 
 \left(\int_K |\langle y,u \rangle|^p dy\right)^\frac{n}{p} \geq  \\
 && 
 \left( \frac{2\ |B^{n-1}_2| \left[2 (1-\varepsilon) \ r \ h_K(u) \right]^\frac{n-1}{2}}
 {\left[ 1-\frac{\varepsilon\left(h_K(u)+1 -\varepsilon\right)}{2 (1-\varepsilon) r}\right]^\frac{n-1}{2}}  \right)^\frac{n}{p}\  
 \left( \int_{(1-\varepsilon) (h_K(u)- \Delta_\varepsilon)}^{h_K(u)}  t^p \left(1-\frac{t}{h_K(u)}\right)^\frac{n-1}{2} 
dt \right)^\frac{n}{p}   = \\
&& \left( \frac{ |B^{n-1}_2| \left((1-\varepsilon) \ r \right)^\frac{n-1}{2}\ \left[2 h_K(u) \right]^\frac{n+1}{2}}
 {\left[ 1-\frac{\varepsilon\left(h_K(u)+1 -\varepsilon\right)}{2 (1-\varepsilon) r}\right]^\frac{n-1}{2}}  \right)^\frac{n}{p} h_K(u)^n 
 \left( \int_{(1-\varepsilon) (1- \frac{\Delta_\varepsilon}{h_K(u)})}^{1}  v^p (1-v)^\frac{n-1}{2} 
dv \right)^\frac{n}{p}
\end{eqnarray*}
Now we apply Lemma \ref{klar} to the function $f(v)=(1-v)^\frac{n-1}{2}$. 
$f$ is $C^2$ and $v_p=\frac{1}{1+\frac{n-1}{2p}}$. Thus  Lemma \ref{klar} holds. 
$v_p$  of Lemma \ref{klar} is an increasing function of $p$ and $\lim_{p\rightarrow \infty} v_p =1$. 
Hence, for $\varepsilon >0$ given there exists $p_\varepsilon=p_{\varepsilon, \Delta_\varepsilon}$
namely $p_\varepsilon \geq \frac{(n-1)\left(h_K(u)  - \Delta_\varepsilon\right)}{2 \Delta _\varepsilon}$,  such that
for all $p \geq p_\varepsilon$, $v_p \geq \frac{h_K(u) - \Delta_\varepsilon}{h_K(u)}$. In addition, we also choose $p_\varepsilon$ so large so  that $p_\varepsilon \geq \frac{1}{\varepsilon^3}$.
Thus
\begin{eqnarray*}
\frac{h_{Z_p(K)}^{n}(u)}{h_{K}^{n}(u)} &\geq &
\left( \frac{ |B^{n-1}_2| \left((1-\varepsilon) \ r\right)^\frac{n-1}{2}  \left[2h_K(u) \right]^\frac{n+1}{2}}{\left(1+C e^{-\frac{c }{\varepsilon}}\right)\  \left[ 1-\frac{\varepsilon\left(h_K(u)+1 -\varepsilon\right)}{2 (1-\varepsilon) r}\right]^\frac{n-1}{2}} \right)^\frac{n}{p} 
 \left( \int_{0}^{1} v^p (1-v)^\frac{n-1}{2} dv 
\right)^\frac{n}{p}.
\end{eqnarray*}
Now
\begin{eqnarray} \label{equation1}
&& \hskip -10mm \left( \frac{ |B^{n-1}_2| \left((1-\varepsilon) \ r\right)^\frac{n-1}{2}  \left[2h_K(u) \right]^\frac{n+1}{2}}{\left(1+C e^{-\frac{c }{\varepsilon}}\right)\  \left[ 1-\frac{\varepsilon\left(h_K(u)+1 -\varepsilon\right)}{2 (1-\varepsilon) r}\right]^\frac{n-1}{2}} \right)^\frac{n}{p}  = 
1 +\frac{n}{p} \log\left(\frac{|B^{n-1}_2| \left((1-\varepsilon)\  r\right)^\frac{n-1}{2}  \left[2h_K(u) \right]^\frac{n+1}{2}}{\left(1+ C e^{-\frac{c }{\varepsilon}}\right)\  \left[ 1-\frac{\varepsilon\left(h_K(u)+1 -\varepsilon\right)}{2 (1-\varepsilon) r}\right]^\frac{n-1}{2} } \right)
 \nonumber \\
&&\hskip -5mm +  \frac{1}{2} \left( \frac{n}{p} \log\left(\frac{|B^{n-1}_2| \left((1-\varepsilon)\  r\right)^\frac{n-1}{2}  \left[2h_K(u) \right]^\frac{n+1}{2}}{\left(1+ C e^{-\frac{c }{\varepsilon}}\right)\ \left[ 1-\frac{\varepsilon\left(h_K(u)+1 -\varepsilon\right)}{2 (1-\varepsilon) r}\right]^\frac{n-1}{2}  } \right) \right)^2 \pm o(p^2).
\end{eqnarray}
Together with Lemma \ref{lemma1} (ii) (for $a=0$) we then get: For $\varepsilon >0$ given, there exists $p_\varepsilon$  such that for all $p \geq p_\varepsilon$

\begin{eqnarray}
&&\frac{h_{Z_p(K)}^{n}(u)}{h_{K}^{n}(u)} \geq \nonumber \\
&&1 
- \frac{n(n+1)}{2p} \log p
+ \frac{n}{2p} \log\left(\frac{\pi^{n-1} \left((1-\varepsilon) r\right)^{n-1} \left[2h_K(u) \right]^{n+1}}{\left(1+ C e^{-\frac{c }{\varepsilon}} \right)^2 \ \left[ 1-\frac{\varepsilon\left(h_K(u)+1 -\varepsilon\right)}{2 (1-\varepsilon) r}\right]^{n-1}\bigg)} \right) +\nonumber \\
&&  \frac{n^2(n+1)^2}{8p^2} (\log p)^2 - \frac{n^2(n+1)}{2p^2} \log\left(
\frac{\pi^{n-1} \left((1-\varepsilon) r\right)^{n-1} \left[2h_K(u) \right]^{n+1}}
{\left(1+ C e^{-\frac{c }{\varepsilon}} \right)^2\ \left[ 1-\frac{\varepsilon\left(h_K(u)+1 -\varepsilon\right)}{2 (1-\varepsilon) r}\right]^{n-1}}  \right) \  \log p \nonumber \\
&& - \  \frac{n(n+1)}{2p^2}\left[  \frac{\left(n^2+3n+6\right)}{4} 
\right]  + \nonumber\\
&& \frac{n^2}{2p^2}\left[ \left(\log \left(\Gamma(\frac{n+1}{2})\right)\right)^2+ 2 \log\left(\frac{\pi^{n-1} \left((1-\varepsilon) r\right)^{n-1} \left[2h_K(u) \right]^{n+1}}
{\left(1+ C e^{-\frac{c }{\varepsilon}} \right)^2 \ \left[ 1-\frac{\varepsilon\left(h_K(u)+1 -\varepsilon\right)}{2 (1-\varepsilon) r}\right]^{n-1}}  \right)  \right]  + \nonumber \\
&& \frac{n^2}{2p^2}\left[    \left(\log\left(   \frac{ |B^{n-1}_2| \left((1-\varepsilon) \ r \right)^\frac{n-1}{2}  \left[2h_K(u) \right]^\frac{n+1}{2}}{\left(1+C e^{-\frac{c }{\varepsilon}}\right )\ \left[ 1-\frac{\varepsilon\left(h_K(u)+1 -\varepsilon\right)}{2 (1-\varepsilon) r}\right]^\frac{n-1}{2}}  \right) \right)^2 \right] \pm o(p^2). \label{h-below}
\end{eqnarray}

Thus
\begin{eqnarray}
&&\frac{p}{\log p}\ \left(1- \frac{h_{Z_p(K)}^{n}(u)}{h_{K}^{n}(u)}\right) \leq \nonumber\\
&& \frac{n(n+1)}{2} - 
\frac{n}{2\log p} \log\left(\frac{\pi^{n-1} \left((1-\varepsilon) r\right)^{n-1} \left[2h_K(u) \right]^{n+1}
} {\left(1+ C e^{-\frac{c }{\varepsilon}} \right)^2  \ \left[ 1-\frac{\varepsilon\left(h_K(u)+1 -\varepsilon\right)}{2 (1-\varepsilon) r}\right]^{n-1}}\right)   \pm o(p). \label{above}
\end{eqnarray}

\vskip 4mm
\noindent
On the other hand, by Lemma \ref{derivative},  the function $f(t)= |K \cap (u^{\perp} +t u)|$ satisfies the assumptions of  Lemma \ref{klar} and $t_p$ is well defined.  Also,  $t_p$ is an increasing function of $p$ and by Lemma \ref{limit=h}, $\lim_{p\rightarrow \infty} t_p = h_K(u)$. 
Hence, for $\varepsilon >0$ given there exists $p_\varepsilon=p_{\varepsilon, \Delta_\varepsilon}$ such that
for all $p \geq p_\varepsilon$, $t_p \geq h_K(u) - \Delta_\varepsilon$. In addition, we also choose $p_\varepsilon$ so large so  that $p_\varepsilon \geq \frac{1}{\varepsilon^3}$. Thus
\begin{eqnarray*}
h_{Z_p(K)}^{n}(u) &= & \left(2 \int_0^{h_K(u)}  t^p |\{y \in K: \langle u,y \rangle =t\}| dt \right)^\frac{n}{p}\\
&\leq &
\left(2 \left(1+ C e^{-c \varepsilon^2 p} \right) \int_{t_p (1-\varepsilon)}^{h_K(u)}  t^p |\{y \in K: \langle u,y \rangle =t\}| dt \right)^\frac{n}{p}\\
& \leq& \left(2 \left(1+ C e^{-\frac{c}{ \varepsilon}} \right) \int_{(1-\varepsilon) (h_K(u) -\Delta_\varepsilon)}^{h_K(u)}  t^p |\{y \in K: \langle u,y \rangle =t\}| dt \right)^\frac{n}{p}
\end{eqnarray*}
$$
\leq
\left(2 \left(1+ C e^{-\frac{c}{ \varepsilon}} \right) \int_{(1-\varepsilon) (h_K(u) -\Delta_\varepsilon)}^{h_K(u)}  t^p | 
\{y \in B \big(x-(1+\varepsilon) r\ u, (1+\varepsilon) r \big)  : \langle u, y \rangle =t \}| dt \right)^\frac{n}{p}.
$$
In the last inequality we have used Lemma \ref{Dupin}.
The latter is
$$
\leq 
\left(2 \left(1+ C e^{-\frac{c}{ \varepsilon}} \right) \int_{0}^{h_K(u)}  t^p | 
\{y \in B \big(x-(1+\varepsilon) r\ u, (1+\varepsilon) r \big)  : \langle u, y \rangle =t \}| dt \right)^\frac{n}{p}.
$$

As above,   we notice  that 
$B(x-(1+\varepsilon) r\ u, (1+\varepsilon) r) \cap \{y \in \mathbb{R}^n: \langle u,y \rangle =t\}$ is a 
$(n-1)$-dimensional Euclidean ball  with radius
$$
\bigg( 2 (1+\varepsilon) r \big(h_K(u)-t\big) \left[ 1-\frac{h_K(u)-t}{2 (1+\varepsilon) r}\right]\bigg)^\frac{1}{2}
$$
which is smaller than or equal
$$
\bigg( 2 (1+\varepsilon) r \big(h_K(u)-t\big) \bigg)^\frac{1}{2}
$$

We continue similar to above and get that 
there exists (a new) $p_\varepsilon$  (chosen larger than the ones previously chosen and larger than $\frac{1}{\varepsilon^3}$) such that for all $p \geq p_\varepsilon$
\begin{eqnarray}
&&\frac{h_{Z_p(K)}^{n}(u)}{h_{K}^{n}(u)} \leq 
1 
- \frac{n(n+1)}{2p} \log p
+ \frac{n}{2p} \log\left(\frac{\pi^{n-1} \left((1+\varepsilon) r\right)^{n-1} \left[2h_K(u) \right]^{n+1}}{(1+ C e^{-\frac{c }{\varepsilon}} )^{-2}} \right) +\nonumber \\
&&  \frac{n^2(n+1)^2}{8p^2} (\log p)^2 - \frac{n^2(n+1)}{2p^2} \log\left(
\frac{\pi^{n-1} \left((1+\varepsilon) r\right)^{n-1} \left[2h_K(u) \right]^{n+1}}
{(1+ C e^{-\frac{c }{\varepsilon}} )^{-2}}  \right) \  \log p \nonumber \\
&& - \  \frac{n(n+1)}{2p^2}\left[  \frac{\left(n^2+3n+6\right)}{4} 
\right]  + \nonumber\\
&& \frac{n^2}{2p^2}\left[ \left(\log \left(\Gamma(\frac{n+1}{2})\right)\right)^2+ 2 \log\left(\frac{\pi^{n-1} \left((1+\varepsilon) r\right)^{n-1} \left[2h_K(u) \right]^{n+1}}
{(1+ C e^{-\frac{c }{\varepsilon}} )^{-2}}  \right)  \right]  + \nonumber \\
&& \frac{n^2}{2p^2}\left[    \left(\log\left(   \frac{ |B^{n-1}_2| \left((1+\varepsilon) \ r \right)^\frac{n-1}{2}  \left[2h_K(u) \right]^\frac{n+1}{2}}{(1+C e^{-\frac{c }{\varepsilon}})^{-1}} \right) \right) \right] \pm o(p^2). \label{h-above}
\end{eqnarray}
Thus
\begin{eqnarray}
&&\frac{p}{\log p}\ \left(1- \frac{h_{Z_p(K)}^{n}(u)}{h_{K}^{n}(u)}\right) \geq \nonumber \\
&& \frac{n(n+1)}{2} - 
\frac{n}{2\log p} \log\left(\frac{\pi^{n-1} \left((1+\varepsilon) r\right)^{n-1} \left[2h_K(u) \right]^{n+1}
} {(1+ C e^{-\frac{c }{\varepsilon}} )^{-2}}\right)   \pm o(p). \label{below}
\end{eqnarray}

\noindent
(\ref{above}) and (\ref{below}) give that
$$ 
\lim_{p\rightarrow \infty} \frac{p}{\log{p}} \left( 1- \frac{h_{Z_p(K)}^{n}(u)}{h_{K}(u)^{n}} \right) = \frac{n(n+1)}{2}.
$$
\noindent
Hence, also  using that, since $|K|=1$, $ h_{Z_p(K)}(u) \rightarrow h_{K}(u)$, 
\begin{eqnarray*}
\lim_{p\rightarrow \infty} \frac{p}{\log{p}} \left(|Z_{p}^\circ(K) |-|K^\circ |\right )&= &\frac{1}{n}\  \int_{S^{n-1}}
\lim_{p\rightarrow \infty} \frac{p}{\log p}\  \frac{1}{h_{Z_p(K)}^{n}(u)}\ \left(1- \frac{h_{Z_p(K)}^{n}(u)}{h_{K}^{n}(u)}\right)d \sigma(u)\\
&=& \frac{1}{n}\  \int_{S^{n-1}} \lim_{p\rightarrow \infty} \frac{1}{h_{Z_p(K)}^{n}(u)}
\lim_{p\rightarrow \infty} \frac{p}{\log p}\   \left(1- \frac{h_{Z_p(K)}^{n}(u)}{h_{K}^{n}(u)}\right)d \sigma(u)\\
&=& \frac{n+1}{2} \int_{S^{n-1}} \frac{1}{h_{K}^{n}(u)} d \sigma(u) \\
&=&  \frac{n(n+1)}{2} |K^\circ |.
\end{eqnarray*}

This finishes (i).
\vskip 3mm
\noindent
(ii)
\begin{eqnarray*}
&&  |Z_{p}^\circ(K) |-|K^\circ | - \frac{n(n+1)\log{p}}{2p} |K^\circ| = \\
&& \frac{1}{n}\  \int_{S^{n-1}} \left(\frac{1}{h_{Z_p(K)}^{n}(u)} - \frac{1}{h_{K}^{n}(u)} - \frac{n(n+1)}{2}\  
\frac{\mbox{log}(p)}{p}\  \frac{1}{h_{K}^n(u)}\right) d \sigma(u)= \\
&& \frac{1}{n}\  \int_{S^{n-1}} \frac{1}{h_{Z_p(K)}^{n}(u)} \left(1 - \frac{h_{Z_p(K)}^{n}(u)}{h_{K}^{n}(u)} - 
\frac{n(n+1)}{2}\  
\frac{\mbox{log}(p)}{p}\  \frac{h_{Z_p(K)}^n(u)}{h_{K}^n(u)}\right) d \sigma(u).
\end{eqnarray*}
Hence
\begin{eqnarray*}
&& \lim _{p \rightarrow \infty} p \left(   |Z_{p}^\circ(K) |-|K^\circ | - \frac{n(n+1)\log{p}}{2p} |K^\circ|   \right)= \\
&& \frac{1}{n}  \int_{S^{n-1}} \lim _{p \rightarrow \infty}   \frac{p}{h_{Z_p(K)}^{n}(u)} \left(1 - \frac{h_{Z_p(K)}^{n}(u)}{h_{K}^{n}(u)} - 
\frac{n(n+1)}{2}\  
\frac{\mbox{log}(p)}{p}\  \frac{h_{Z_p(K)}^n(u)}{h_{K}^n(u)}\right) d \sigma(u),
\end{eqnarray*}
where we have used Lemma \ref{interchange} (ii) and Lebegue's theorem to interchange integration and limit.
By  (\ref{h-below})  we have for all $p \geq p_\varepsilon$ 
\begin{eqnarray*}
&&\left(1 - \frac{h_{Z_p(K)}^{n}(u)}{h_{K}^{n}(u)} - 
\frac{n(n+1)}{2}\  
\frac{\mbox{log}(p)}{p}\  \frac{h_{Z_p(K)}^n(u)}{h_{K}^n(u)}\right)  \leq\\
&&-\frac{n}{2p} \log\left(\frac{\pi^{n-1} \left((1-\varepsilon) r\right)^{n-1} \left[2h_K(u) \right]^{n+1}}{\left(1+ C e^{-\frac{c }{\varepsilon}} \right)^2 \left[ 1-\frac{\varepsilon\left(h_K(u)+1 -\varepsilon\right)}{2 (1-\varepsilon) r}\right]^{n-1}} \right) +
  \frac{n^2(n+1)^2}{8p^2} (\log p)^2 \nonumber \\
&& + \  \frac{n(n+1)}{2p^2}\left[  \frac{\left(n^2+3n+6\right)}{4} 
\right]  - \nonumber\\
&& \frac{n^2}{2p^2}\left[ \left(\log \left(\Gamma(\frac{n+1}{2})\right)\right)^2+ 2 \log\left(\frac{\pi^{n-1} \left((1-\varepsilon) r\right)^{n-1} \left[2h_K(u) \right]^{n+1}}
{\left(1+ C e^{-\frac{c }{\varepsilon}} \right)^2    \left[ 1-\frac{\varepsilon\left(h_K(u)+1 -\varepsilon\right)}{2 (1-\varepsilon) r}\right]^{n-1} }  \right)  \right]  - \nonumber \\
&& \frac{n^2}{2p^2}\left[    \left(\log\left(   \frac{ |B^{n-1}_2| \left((1-\varepsilon) \ r \right)^\frac{n-1}{2}  \left[2h_K(u) \right]^\frac{n+1}{2}}{\left(1+C e^{-\frac{c }{\varepsilon}}\right)  \left[ 1-\frac{\varepsilon\left(h_K(u)+1 -\varepsilon\right)}{2 (1-\varepsilon) r}\right]^\frac{n-1}{2} }  \right) \right)^2 \right] \pm o(p^2)
\end{eqnarray*}
Thus
\begin{eqnarray}
&&p\ \left(1 - \frac{h_{Z_p(K)}^{n}(u)}{h_{K}^{n}(u)} - 
\frac{n(n+1)}{2}\  
\frac{\mbox{log}(p)}{p}\  \frac{h_{Z_p(K)}^n(u)}{h_{K}^n(u)}\right)  \nonumber \leq \\
&&-\frac{n}{2} \log\left(\frac{\pi^{n-1} \left((1-\varepsilon) r\right)^{n-1} \left[2h_K(u) \right]^{n+1}}{\left(1+ C e^{-\frac{c }{\varepsilon}} \right)^2  \left[ 1-\frac{\varepsilon\left(h_K(u)+1 -\varepsilon\right)}{2 (1-\varepsilon) r}\right]^{n-1}} \right) +
  \frac{n^2(n+1)^2}{8p} (\log p)^2 \nonumber \\
&& + \  \frac{n(n+1)}{2p}\left[  \frac{\left(n^2+3n+6\right)}{4} 
\right]  - \nonumber\\
&& \frac{n^2}{2p}\left[ \left(\log \left(\Gamma(\frac{n+1}{2})\right)\right)^2+ 2 \log\left(\frac{\pi^{n-1} \left((1-\varepsilon) r\right)^{n-1} \left[2h_K(u) \right]^{n+1}}
{\left(1+ C e^{-\frac{c }{\varepsilon}} \right)^2  \left[ 1-\frac{\varepsilon\left(h_K(u)+1 -\varepsilon\right)}{2 (1-\varepsilon) r}\right]^{n-1}}  \right)  \right]  - \nonumber \\
&& \frac{n^2}{2p}\left[    \left(\log\left(   \frac{ |B^{n-1}_2| \left((1-\varepsilon) \ r \right)^\frac{n-1}{2}  \left[2h_K(u) \right]^\frac{n+1}{2}}{\left(1+C e^{-\frac{c }{\varepsilon}} \right)  \left[ 1-\frac{\varepsilon\left(h_K(u)+1 -\varepsilon\right)}{2 (1-\varepsilon) r}\right]^\frac{n-1}{2}}  \right) \right)^2 \right] \pm o(p) \label{above2}
\end{eqnarray}
Similarly, using (\ref{h-above}), we get  for all $p \geq p_\varepsilon$ 
\begin{eqnarray}
&&p\ \left(1 - \frac{h_{Z_p(K)}^{n}(u)}{h_{K}^{n}(u)} - 
\frac{n(n+1)}{2}\  
\frac{\mbox{log}(p)}{p}\  \frac{h_{Z_p(K)}^n(u)}{h_{K}^n(u)}\right)  \nonumber \geq \\
&&-\frac{n}{2} \log\left(\frac{\pi^{n-1} \left((1+\varepsilon) r\right)^{n-1} \left[2h_K(u) \right]^{n+1}}{(1+ C e^{-\frac{c }{\varepsilon}} )^{-2}} \right) +
  \frac{n^2(n+1)^2}{8p} (\log p)^2 \nonumber \\
&& + \  \frac{n(n+1)}{2p}\left[  \frac{\left(n^2+3n+6\right)}{4} 
\right]  - \nonumber\\
&& \frac{n^2}{2p}\left[ \left(\log \left(\Gamma(\frac{n+1}{2})\right)\right)^2+ 2 \log\left(\frac{\pi^{n-1} \left((1+\varepsilon) r\right)^{n-1} \left[2h_K(u) \right]^{n+1}}
{(1+ C e^{-\frac{c }{\varepsilon}} )^{-2}}  \right)  \right]  - \nonumber \\
&& \frac{n^2}{2p}\left[    \left(\log\left(   \frac{ |B^{n-1}_2| \left((1+\varepsilon) \ r \right)^\frac{n-1}{2}  \left[2h_K(u) \right]^\frac{n+1}{2}}{(1+C e^{-\frac{c }{\varepsilon}})^{-1}}  \right) \right)^2 \right] \pm o(p) \label{below2}
\end{eqnarray}
\par
\noindent
(\ref{above2}) and (\ref{below2}) give that
\begin{eqnarray*}
&&\lim_{p\rightarrow \infty} p\ \left(1 - \frac{h_{Z_p(K)}^{n}(u)}{h_{K}^{n}(u)} - 
\frac{n(n+1)}{2}\  
\frac{\mbox{log}(p)}{p}\  \frac{h_{Z_p(K)}^n(u)}{h_{K}^n(u)}\right) \\
&&= 
- \frac{n}{2} \log\left(\pi^{n-1}  r^{n-1} \left[2h_K(u) \right]^{n+1} \right).
\end{eqnarray*} 
\vskip 3mm
The  limit  $\ \lim _{p \rightarrow \infty}  p \left( |Z_p^\circ(K)| - |K^\circ|- \frac{n(n+1)}{2p} \log p \ |Z_p^\circ(K)| \right)$ is computed similarly.

\vskip 5mm

\section{Applications}
 
\medskip
The fact that  $\Omega_K$ can  be expressed in different ways allows us to compute the integral in the 
next proposition.
\vskip 3mm

\begin{proposition}
Let $1 <  r <\infty$ and let $B^n_r$ be the $l_r^n$- unit ball and let $(B_{r}^{n-1})^{+}$ be the set of all
vectors in $B_{r}^{n-1}$ having nonnegative
coordinates.
 Then
\begin{eqnarray*}
&&\int_{(B_{r}^{n-1})^{+}}
\prod_{i=1}^{n-1}|x_{i}|^{r-2}
 \ \log\left[ (r-1)^{n-1}\prod_{i=1}^{n}|x_{i}|^{r-2}\right] 
x_{n}^{-1}
dx_1\dots dx_{n-1} =  \\
&&\frac{n}{r^{n-1}} \frac{\big(\Gamma(\frac{r-1}{r})\big)^n}{\Gamma(\frac{n(r-1)}{r})} \left[   \frac{n(r-2)}{r} \left( \frac{\Gamma^{\prime}(\frac{r-1}{r})}{\Gamma(\frac{r-1}{r})} - \frac{\Gamma^{\prime}(n\frac{r-1}{r})}{\Gamma(n\frac{r-1}{r})}\right)\bigg) + (n-1) \ \log r  \right]
\end{eqnarray*}

\end{proposition}
\vskip 3mm
\noindent
{\bf Proof.}
\par
\noindent
In Chapter 3 it was shown that 
\begin{eqnarray*}
 \log{\Omega_{K}} = 
 - \frac{n}{as_{\infty}(K)} \int_{\partial K} \frac{\kappa_{K}(x)}{\langle x, N_{K}(x) \rangle^{n}}  \log{ \frac{\kappa_{K}(x)}{\langle x, N_{K}(x) \rangle^{n+1} } } d\mu_{K}(x).
\end{eqnarray*}

\noindent
We apply this formula to  $K=B^n_r$, 
$1 <  r <\infty$. It was also shown in Chapter 3 that
\begin{equation*}
\log{\Omega_{B^n_r}} = - n \left[   \frac{n(r-2)}{r} \left( \frac{\Gamma^{\prime}(\frac{r-1}{r})}{\Gamma(\frac{r-1}{r})} - \frac{\Gamma^{\prime}(n\frac{r-1}{r})}{\Gamma(n\frac{r-1}{r})}\right)\bigg) + (n-1) \ \log r  \right]
\end{equation*}
\noindent
The curvature  at a boundary point of $B^n_r$ is (see \cite{SW5})
$$
\kappa(x)=\frac{(r-1)^{n-1}\prod_{i=1}^{n}|x_{i}|^{r-2}}
{(\sum_{i=1}^{n}|x_{i}|^{2r-2})^{\frac{n+1}{2}}}
$$
and the normal is (see \cite{SW5})
$$
N_{\partial B^n_r}(x)=\frac{(\mbox{sgn}(x_{1})|x_{1}|^{r-1},\dots,
\mbox{sgn}(x_{n})|x_{n}|^{r-1})}{(\sum_{i=1}^{n}|x_{i}|^{2r-2}
)^{\frac{1}{2}}}.
$$
Thus we get - where $B^n_{r^\prime} $ is the polar of $B^n_{r}$, i.e. $r^\prime $ is the conjugate exponent of $r$ - 
\begin{eqnarray*}
&& n \left[   \frac{n(r-2)}{r} \left( \frac{\Gamma^{\prime}(\frac{r-1}{r})}{\Gamma(\frac{r-1}{r})} - \frac{\Gamma^{\prime}(n\frac{r-1}{r})}{\Gamma(n\frac{r-1}{r})}\right)\bigg) + (n-1) \ \log r  \right] |B^n_{r^\prime}| = \\
&& \int_{\partial B_{r}^{n}}
\frac{((r-1)^{n-1}\prod_{i=1}^{n}|x_{i}|^{r-2}}
{(\sum_{i=1}^{n}|x_{i}|^{2r-2})^{\frac{1}{2}}} \ \log\left[ (r-1)^{n-1}\prod_{i=1}^{n}|x_{i}|^{r-2}\right] 
d\mu_{\partial B_{r}^{n}}(x).
\end{eqnarray*}
  Now we integrate with respect to the variables $x_{1},\dots,x_{n-1}$.
The volume of a surface element in the plane of the first $n-1$
coordinates equals the volume of the corresponding surface element
on $\partial B_{r}^{n}$ times
$$
|<e_{n},N_{\partial B_{r}^{n}}(x)>|
=\frac{|x_{n}|^{r-1} }{(\sum_{i=1}^{n}|x_{i}|^{2r-2})^{\frac{1}{2}}}.
$$
Thus, with $(B_{r}^{n-1})^{+}$ being the set of all
vectors in $B_{r}^{n-1}$ having nonnegative
coordinates, 
\begin{eqnarray*}
 2^{n}(r-1)^{n-1} \int_{(B_{r}^{n-1})^{+}}
\prod_{i=1}^{n}|x_{i}|^{r-2}
 \ \log\left[ (r-1)^{n-1}\prod_{i=1}^{n}|x_{i}|^{r-2}\right] 
x_{n}^{1-r}
dx_1\dots dx_{n-1}   \nonumber \\
= 2^{n}(r-1)^{n-1} \int_{(B_{r}^{n-1})^{+}}
\prod_{i=1}^{n-1}|x_{i}|^{r-2}
 \ \log\left[ (r-1)^{n-1}\prod_{i=1}^{n}|x_{i}|^{r-2}\right] 
x_{n}^{-1}
dx_1\dots dx_{n-1}   \nonumber \\
=  2^{n}(r-1)^{n-1} \frac{n}{r^{n-1}} \frac{\big(\Gamma(\frac{r-1}{r})\big)^n}{\Gamma(\frac{n(r-1)}{r})} \left[   \frac{n(r-2)}{r} \left( \frac{\Gamma^{\prime}(\frac{r-1}{r})}{\Gamma(\frac{r-1}{r})} - \frac{\Gamma^{\prime}(n\frac{r-1}{r})}{\Gamma(n\frac{r-1}{r})}\right)\bigg) + (n-1) \ \log r  \right],
\end{eqnarray*}
where we have also used that
$$
|B^n_{r^\prime}| =\frac{2^{n}(r-1)^{n-1} }{n \  r^{n-1}} \frac{\big(\Gamma(\frac{r-1}{r})\big)^n}{\Gamma(\frac{n(r-1)}{r})}.
$$

\vskip 3mm
\noindent
There are still other  ways how $\Omega_K$ can be expressed. Similar to Theorem \ref{theorem1}, $\Omega_K$
appears in the asymptotic behavior of the volume of certain  surface bodies and illumination surface bodies \cite{WY1}. We show the result
for the surface bodies. For the illumination surface bodies it is done similarly.
\par
The surface bodies, a variant of the floating bodies,  were introduced 
 in \cite{SW4, SW5} as follows
\vskip 3mm
\noindent
{\bf Definition}
\par
Let $s\geq 0$ and $f:\partial K\rightarrow \mathbb{R}$ be a
nonnegative, integrable function. The {\it surface body $K_{f,s}$}
is the intersection of all the closed half-spaces $H^+$ whose
defining hyperplanes $H$ cut off a set of $f \mu_K$-measure less
than or equal to $s$ from $\partial K$. More precisely,
\begin{equation*}\label{def1}
K_{f,s}= \bigcap _{\int_{\partial K \cap H^-}f d\mu_K \leq s} H^+.
\end{equation*}

\vskip 3mm

\noindent
\begin{proposition}
Let $K$ be a symmetric convex body in $\mathbb{R}^n$ that is in $C^2_+$. 
 
\begin{eqnarray*}
 d_{n}  \mbox{lim}_{s \rightarrow 0} \frac{|K| - |K_{f,s}|}{s^{\frac{2}{n-1}}} &=& \\
&& \hskip -10mm  \int_{\partial K} \frac{ \kappa(x)}{\langle x, N(x) \rangle^n} \   \log \left(\frac{\kappa(x)}{
\langle x, N(x) \rangle^{n+1}} \right) d\mu(x) =   |K^\circ|\  \log\frac{1}{\Omega_{K}}. 
\end{eqnarray*}
\vskip 2mm
\noindent
where $K_{f,s}$ is the surface body of $K$  for the function 
$$
f=\frac{\langle x,N_K(x)\rangle^{\frac{n(n-1)}{2}}}{\kappa^\frac{n-2}{2}} \ \bigg(\log \  \big(\frac {\kappa}{ \langle x,N_K(x)\rangle^{n+1}}\big)\bigg)^{-\frac{n-1}{2}}
$$ and  where $d_n=2 \bigg(|
B^{n-1}_2|\bigg)^\frac{2}{n-1}$.
\end{proposition}
\vskip 3mm
\noindent
{\bf Proof.} 
\par
\noindent
The proof follows immediately from the following  formula which was proved in \cite{SW5} (Theorem 14)
$$
d_{n} \lim_{s \to 0}
\frac{|K|-| K_{f,s}  |}
{s^\frac{2}{n-1}}=
\int_{\partial K} \frac{\kappa^\frac{1}{n-1}}
{f^\frac{2}{n-1}}d\mu_{\partial K}.
$$

\vskip 5mm

\section{Appendix: Calculations with $\Gamma$-functions.}

\vskip 4mm 
\noindent 
For  $x,y>0$,  $\Gamma(x):= \int_{0}^{\infty} \lambda^{x-1} e^{-\lambda} d\lambda $ is the Gamma function and 
$ B(x,y):= \int_{0}^{1} \lambda^{x-1} (1-\lambda)^{y-1} d \lambda = \frac{\Gamma(x) \Gamma(y)}{\Gamma(x+y)} $ is the Beta function.
\par
Recall that 
we write $f(p)=g(p)\pm o(p)$, if there exists a function $h(p)$ such that $f(p)=g(p)+ h(p)$ and $\lim_{p\rightarrow \infty} ph(p) =0$
and similarly,   $f(p)=g(p)\pm o(p^2)$, if   there exists a function $h(p)$ such that $f(p)=g(p)+ h(p)$ and $\lim_{p \rightarrow \infty} p^2h(p) =0$

\vskip 3mm
\noindent
We will  frequently use:
For $x \rightarrow \infty$, 
\begin{eqnarray}\label{gamma}
\Gamma(x)= \sqrt{2\pi} \ x^{x-\frac{1}{2}}\  e^{-x}\ \left[1+ \frac{1}{12x}  + \frac{1}{288x^2}  \pm o(x^2)\right].
\end{eqnarray}
For every $z,w>0$ 
$$ 
z^{1/p} = 1 + \frac{\log{z}}{p} + \frac{(\mbox{log}z)^2}{2 p^2} \pm o(p^2) 
$$
and
$$ 
(p+z)^{w/p} = 1 + \frac{w}{p}\log{p} + \frac{w^2 (\mbox{log}z)^2}{2 p^2} +\frac{wz}{p^2} \pm o(p^2).
$$

\noindent 
Note  that if $f(p)^2 = o(p)$ then $(1+f(p))(1-f(p)) =1 \pm o(p)$, which means that
$$ \frac{1}{1+f(p)}= 1-f(p) \pm o(p). $$

\noindent 
Also 

$$ \frac{a}{p+b}= \frac{a}{p} -\frac{ab}{p^2} \pm o(p^2). $$
\vskip 3mm
\noindent
{\bf Proof of Lemma \ref{lemma1}}
\vskip 3mm
\noindent
(i)
We use 
(\ref{gamma}) and get 
\begin{eqnarray*}
&&\left(B\big(p+1,\frac{n+1}{2}\big) \right)^{\frac{n}{p}} 
=
 \left( \frac{\Gamma(p+1)}{\Gamma(p+1+\frac{n+1}{2})} \Gamma(\frac{n+1}{2})\right)^{\frac{n}{p}}\\
 &&= 
 \left(  \frac{ \Gamma(\frac{n+1}{2})\  e^{\frac{n+1}{2}} \ (p+1)^{p+\frac{1}{2}} \left[1+ \frac{1}{12(p+1)}  + \frac{1}{288(p+1)^2}  \pm o(p^2)\right] }{ (p+1+\frac{n+1}{2})^{p+1+\frac{n}{2}} \left[1+ \frac{1}{12(p+1+\frac{n+1}{2})} +  \frac{1}{288(p+1+\frac{n+1}{2})^2}  \pm o(p^2)\right]} \right)^{\frac{n}{p}} \\
&&= \left( \Gamma\left(\frac{n+1}{2}\right)\ e^{\frac{n+1}{2}} \right)^{\frac{n}{p}} 
\left( \frac{p+1}{p+1+\frac{n+1}{2}}\right)^{\frac{n}{p}(p+\frac{1}{2})} 
\left( \frac{1}{p+1+\frac{n+1}{2}}\right)^{\frac{n(n+1)}{2p}} \times \\
&& \hskip 37mm  \left(\frac{1+ \frac{1}{12(p+1)}  + \frac{1}{288(p+1)^2}  \pm o(p^2)} {1+ \frac{1}{12(p+1+\frac{n+1}{2})} +  \frac{1}{288(p+1+\frac{n+1}{2})^2}  \pm o(p^2)}\right)^\frac{n}{p}
 \end{eqnarray*}

\noindent 
Note that 
$$
\left(\frac{1+ \frac{1}{12(p+1)}  + \frac{1}{288(p+1)^2}  \pm o(p^2)} {1+ \frac{1}{12(p+1+\frac{n+1}{2})} +  \frac{1}{288(p+1+\frac{n+1}{2})^2}  \pm o(p^2)}\right)^\frac{n}{p} = 1\pm o(p^2).
$$
Also
\begin{eqnarray*}
&& \left(\Gamma\left(\frac{n+1}{2}\right)\ e^{\frac{n+1}{2}} \right)^\frac{n}{p} \\
&&= 1
+  \frac{n}{p} \left[\frac{n+1}{2}  + \mbox{log} \left(\Gamma(\frac{n+1}{2})\right)\right]  + 
\frac{n^2}{2p^2} \left[\frac{n+1}{2}  + \mbox{log} \left(\Gamma(\frac{n+1}{2})\right)\right] ^2
\pm o(p^2),
\end{eqnarray*} 
\par

\begin{eqnarray*}
\left( \frac{1}{1+\frac{n+1}{2(p+1)}}\right)^{n(1+\frac{1}{2p})} 
&= &
\left( \frac{1}{1+\frac{n+1}{2(p+1)}}\right)^{n}
\ e^{- \frac{n}{2p} \log\left(1+ \frac{n+1}{2p+2}\right)} \\
&= &1- \frac{n(n+1)}{2p}  + \frac{n(3+5n+3n^2+n^3)}{8p^2}\pm o(p^2)
\end{eqnarray*} 
and
\begin{eqnarray*}
&&\left( \frac{1}{p+1+\frac{n+1}{2}}\right)^{\frac{n(n+1)}{2p}} = \ e^{- \frac{n(n+1)}{2p} \log\left(p+ \frac{n+3}{2}\right)} \\
&&=1- \frac{n(n+1)}{2p} \log p + \frac{n^2(n+1)^2}{8p^2} (\log p)^2 - \frac{n(n+1)(n+3)}{4p^2}  \pm o(p^2).
\end{eqnarray*} 
Hence
\begin{eqnarray*} 
&&\left(B\big(p+1,\frac{n+1}{2}\big) \right)^{\frac{n}{p}}   =  \bigg(1\pm o(p^2) \bigg) \\
&& \left(1+  \frac{n}{p} \left[\frac{n+1}{2}  + \mbox{log} \left(\Gamma(\frac{n+1}{2})\right)\right] 
+ 
\frac{n^2}{2p^2} \left[\frac{n+1}{2}  + \mbox{log} \left(\Gamma(\frac{n+1}{2})\right)\right] ^2  \pm o(p^2)\right)\\
&& \left( 1- \frac{n(n+1)}{2p} + \frac{n(3+5n+3n^2+n^3)}{8p^2}  \pm o(p^2)\right) \\
&&\left(1- \frac{n(n+1)}{2p} \log p + \frac{n^2(n+1)^2}{8p^2} (\log p)^2 - \frac{n(n+1)(n+3)}{4p^2}  \pm o(p^2) \right)\\
&&= 1- \frac{n(n+1)}{2p} \log p + \frac{n}{p} \mbox{log} \left(\Gamma(\frac{n+1}{2})\right) +
\frac{n^2(n+1)^2}{8p^2} (\log p )^2 \\
&& - \frac{n^2(n+1)}{2p^2}   \mbox{log}  \left(\Gamma(\frac{n+1}{2})\right) \log p \\
&&+\frac{n}{2p^2}\left[n\left(\mbox{log} \left(\Gamma(\frac{n+1}{2})\right)\right)^2- \frac{n+1}{4}\left(n(n+1)+2(n+3)\right) \right]\pm o(p^2).
\end{eqnarray*} 
\vskip 3mm
\noindent
(ii)
\begin{eqnarray*} 
&& \left( \int_{0}^{1} u^p (1-u)^\frac{n-1}{2} \left( 1-a\left(1-u \right)\right)^\frac{n-1}{2} du 
\right)^\frac{n}{p} \\
&&=
\bigg( \int_{0}^{1} u^p (1-u)^\frac{n-1}{2} \bigg[ 1-
{\frac{n-1}{2} \choose 1}\ 
a\  (1-u) +
{\frac{n-1}{2} \choose 2 }a^2 (1-u)^2 \pm \dots
\bigg] du 
\bigg)^\frac{n}{p} \\
&&= 
\left( B\big(p+1, \frac{n+1}{2}\big) \right)^\frac{n}{p} 
 \bigg[ 1- {\frac{n-1}{2} \choose 1 } \ a \ B_3  \ + {\frac{n-1}{2} \choose 2 } \ a^2 \   B_5 -  {\frac{n-1}{2} \choose 3 } \ a^3  \   B_7 \pm \dots
 \bigg] ^\frac{n}{p} \\
&&= 
\left( B\big(p+1, \frac{n+1}{2}\big) \right)^\frac{n}{p} 
\mbox{exp}\bigg\{\frac{n}{p} \log\bigg[ 1- {\frac{n-1}{2} \choose 1} \ a \ B_3 \  +  {\frac{n-1}{2} \choose 2} \ a^2 \   B_5 
\pm \dots
 \bigg] \bigg\} \\
&&= 
\left( B\big(p+1, \frac{n+1}{2}\big) \right)^\frac{n}{p} \times \\
&&\hskip 12mm \bigg[ 1- \frac{n}{p}  \bigg\{ {\frac{n-1}{2} \choose 1} \ a \ B_3 \  -  {\frac{n-1}{2} \choose 2} \ a^2 \   B_5 +\frac{1}{2}  \left({\frac{n-1}{2} \choose 1}\right)^2  a^2 \   B_3^2 \pm \dots \bigg\} \dots \bigg]
\end{eqnarray*} 
where for $3 \leq k \leq n-2$ and for a constant $c$
$$
B_k=\frac{B\big(p+1, \frac{n+k}{2}\big) }{B\big(p+1, \frac{n+1}{2}\big) }
=\frac{\Gamma(\frac{n+k}{2})}{\Gamma(\frac{n+1}{2})} \frac{1}{p^\frac{k-1}{2}} \left(1+ \frac{c}{p} \\ \pm o(p)\right)
$$
Hence, together with (i), 
\begin{eqnarray*}
&&\left( \int_{0}^{1} u^p (1-u)^\frac{n-1}{2} \left( 1-a\left(1-u \right)\right)^\frac{n-1}{2} du 
\right)^\frac{n}{p}= 
 1- \frac{n(n+1)}{2p} \log p + \\
 &&\frac{n}{p} \mbox{log} \left(\Gamma(\frac{n+1}{2})\right) +
\frac{n^2(n+1)^2}{8p^2} (\log p )^2
- \frac{n^2(n+1)}{2p^2}   \mbox{log}  \left(\Gamma(\frac{n+1}{2})\right) \log p +\\
&&\frac{n}{2p^2}\left[n\left(\mbox{log} \left(\Gamma(\frac{n+1}{2})\right)\right)^2- \frac{(n+1)\left(n^2+3n+6\right)}{4}
-2  {\frac{n-1}{2} \choose 1} \ a \  \frac{\Gamma(\frac{n+3}{2})}{\Gamma(\frac{n+1}{2})}\right] +\\
&&\frac{n}{2p^2}\left[n\left(\mbox{log} \left(\Gamma(\frac{n+1}{2})\right)\right)^2- \frac{(n+1)\left(n^2+3n+6\right)}{4}
-(n+1) {\frac{n-1}{2} \choose 1} \ a \right] \\
&&\pm o(p^2).
\end{eqnarray*}

\vskip 2mm 
\noindent 
Grigoris Paouris\\
{\small Department of Mathematics}\\
{\small Texas A \& M University}\\
{\small College Station, TX, , U. S. A.}\\
{\small \tt }\\ \\
\noindent
\and 
Elisabeth Werner\\
{\small Department of Mathematics \ \ \ \ \ \ \ \ \ \ \ \ \ \ \ \ \ \ \ Universit\'{e} de Lille 1}\\
{\small Case Western Reserve University \ \ \ \ \ \ \ \ \ \ \ \ \ UFR de Math\'{e}matique }\\
{\small Cleveland, Ohio 44106, U. S. A. \ \ \ \ \ \ \ \ \ \ \ \ \ \ \ 59655 Villeneuve d'Ascq, France}\\
{\small \tt elisabeth.werner@case.edu}\\ \\

\end{document}